\documentclass[12pt]{amsart}
\usepackage{amsmath, amsfonts, amssymb}
\newcommand{\dal}{\square}
\newcommand{\R}{{\mathbf R}}
\newcommand{\ve}{\varepsilon}
\newcommand{\pa}{\partial}

\newcommand{\jb}[1]{\left\langle #1 \right\rangle}

\newcommand{\Hi}{{\mathcal H}}
\DeclareMathOperator{\supp}{\rm supp}
\newtheorem{theorem}{Theorem}[section]
\newtheorem{proposition}[theorem]{Proposition}

\newtheorem{lemma}[theorem]{Lemma}

\numberwithin{equation}{section}

\usepackage{color}

\title[Mixed problem for the wave equation]
{The rate of convergence to the asymptotics for the wave
equation in an exterior domain}

\author[S.~Katayama]{Soichiro Katayama}
\thanks{The first author is partially supported by Grant-in-Aid for Scientific Research (C)
(No.~20540211), JSPS}
\address{Department of Mathematics, Wakayama University, 930 Sakaedani, Wakayama 640-8510, Japan}
\email{katayama@center.wakayama-u.ac.jp}
\author[H.~Kubo]{Hideo Kubo}
\address{Division of Mathematics,
Graduate School of Information Sciences,
Tohoku University,
Sendai 980-8579, Japan}
\email{kubo@math.is.tohoku.ac.jp}

\date{}

\begin{document}

\begin{abstract}
In this paper we consider the mixed problem for the wave equation
exterior to a non-trapping obstacle in odd space dimensions.
We derive a rate of the convergence of the solution for the mixed problem
to a solution for the Cauchy problem.
As a by-product, we are able to find out the radiation field of solutions to the mixed problem in terms of the scattering data.
\end{abstract}

\maketitle

%%%%%%%%%%%%%%%%%%%%%%%%%%%%%%%%%%
\section{Introduction}
%%%%%%%%%%%%%%%%%%%%%%%%%%%%%%%%%%

This paper is concerned with the global behavior of solutions to 
the mixed problem for the wave equation in an exterior domain\,:
\begin{align}\label{eq}
& (\partial_t^2-\Delta) u(t,x) =0, & (t,x) \in (0,T)\times \Omega,
\\ \label{dc}
& u(t,x)=0, & (t,x) \in (0,T)\times \partial\Omega,
\\ \label{id}
& u(0,x)=f_0(x),\ (\partial_t u)(0,x)=f_1(x), & x\in \Omega,
\end{align}
where $\Omega={\mathbf R}^n \setminus  \overline{\mathcal O}$, 
and ${\mathcal O}$ is a bounded open set in ${\mathbf R}^n$
with smooth boundary.
Throughout this paper, we suppose that $n$ is an odd integer with $n\ge 3$.
We assume that $\Omega$ is connected and that 
the initial data $\vec{f}=(f_0, f_1)$ belongs to the associated energy space
${\mathcal H}_D(\Omega)$.
Here and in the following, for an open set $Y \subset \R^n$, 
$\Hi_D(Y)$ stands for the completion of $(C^\infty_0(Y))^2$
with respect to $\|\vec{f}\|_{\Hi_D(Y)}=\|\nabla f_0\|_{L^2(Y)}+\|f_1\|_{L^2(Y)}$.
$U(t)$ denotes the propagator of the mixed problem \eqref{eq} with \eqref{dc} and \eqref{id}; in other words, we define
$$
U(t)\vec{f}=\left(u(t, \cdot), \pa_t u (t, \cdot) \right)
$$
for $\vec{f}\in{\mathcal H}_D(\Omega)$,  where $u$ is the solution to \eqref{eq}--\eqref{id}.

It is well known that the asymptotic behavior of the solution to the above problem 
is approximated by a solution to the Cauchy problem.
More precisely, for a given initial data
$\vec{f} \in {\mathcal H}_D(\Omega)$ there exists uniquely a scattering data 
$\vec{f}_+ \in {\mathcal H}_D({\mathbf R}^n)$ such that
\begin{align}\label{ki3}
 \| U(t)\vec{f}-U_0(t)\vec{f}_+ \|_{{\mathcal H}_D(\Omega)} \to 0
\quad (t \to \infty),
\end{align}
where, for $\vec{g}=(g_0, g_1)\in {\mathcal H}_D(\R^n)$,
$U_0(t)\vec{g}$ is given by
$$
U_0(t)\vec{g}=\left(u_0(t, \cdot), \pa_t u_0(t, \cdot)\right)
$$
with $u_0$ being the solution to the Cauchy problem
\begin{align}\label{ki1}
& (\partial_t^2-\Delta) u_0 (t,x) =0, & (t,x) \in (0,T)\times {\mathbf R}^n,
\\ \label{ki2}
& u_0(0,x)=g_0(x),\ (\partial_t u_0)(0,x)=g_1(x), & x\in {\mathbf R}^n.
\end{align}
On the other hand, the local energy of $U(t)\vec{f}$ decays to zero
as $t$ tends to infinity. Namely, for any $R>0$ and any $\vec{f}\in {\mathcal H}_D(\Omega)$, we have
\begin{equation}
\lim_{t\to\infty} \int_{ \{x \in \Omega\,;\,|x|<R\} }
 \left\{|\pa_t u(t,x)|^2+|\nabla u(t,x)|^2\right\} dx=0,
\end{equation}
where $u$ is the solution to \eqref{eq}--\eqref{id}.

In view of these facts, we see that the main part of the perturbed wave $U(t)\vec{f}$ escapes from
any ball with a fixed radius as $t \to \infty$ and that it approaches to some unperturbed wave
$U_0(t)\vec{f}_+$ in the sense of the energy.
However, to our knowledge, the rate of the convergence 
in \eqref{ki3}
is not found explicitly in the literature.
Therefore, it is natural to ask at which rate the perturbed wave tends to an unperturbed wave.
In addition, we are interested in the regularity and decay properties of the scattering data $\vec{f}_+$.
Namely, we wish to know whether the scattering data becomes smoother and decays
faster at the spatial infinity or not, if the initial data is smooth and compactly supported.
This consideration might be useful for the application to the nonlinear 
wave equation in an exterior domain. For instance, we are able to obtain a precise lower bound of the lifespan
in our forthcoming paper.

Here we introduce notation in order to state our main result.
Let $m$ be a nonnegative integer and $Y$ be an open set in $\R^n$.
We set ${\mathcal H}^{m}(Y)=H^{m+1}(Y)\times H^{m}(Y)$
and $\|\vec{f}\|_{\Hi^m(Y)}=\|f_0\|_{H^{m+1}(Y)}+\|f_1\|_{H^m(Y)}$
for $\vec{f}=(f_0, f_1) \in \Hi^m(Y)$.
Similarly, we put ${\mathcal W}^{m,\infty}(Y)=W^{m+1,\infty}(Y)\times W^{m,\infty}(Y)$
and $\|\vec{f}\|_{{\mathcal W}^{m,\infty}(Y)}=
\|f_0\|_{W^{m+1,\infty}(Y)}+\|f_1\|_{W^{m,\infty}(Y)}$
for $\vec{f} \in {\mathcal W}^{m,\infty}(Y)$.
Here $H^m(Y)$ (resp. $W^{m,\infty}(Y))$ stands for the Sobolev space based on 
$L^2(Y)$ (resp. $L^\infty(Y)$).
In addition, we denote by 
$X^m(\Omega)$ the set of all $\vec{f}=(f_0, f_1) \in {\mathcal H}^{m}(\Omega)$
satisfying the compatibility condition of the $m$-th order for the problem \eqref{eq}--\eqref{id},
that is $f_j=0$ {on} $\partial\Omega$ for any $j=0,\dots,m$,
where we have set
\begin{equation}\label{data+}
f_j(x)\equiv \Delta f_{j-2}(x)
\quad \mbox{for \ $x \in \overline{\Omega}$ \ and \ $j\ge 2$}.
\end{equation}
Besides, we put ${\mathcal H}^\infty(Y)=\bigcap_{m=0}^\infty {\mathcal H}^m(Y)$
and $X^\infty(\Omega)=\bigcap_{m=0}^\infty X^m(\Omega)$.

We will use the notation $\pa_j=\pa_{x_j}$ for $1\le j\le n$, and
$\pa_x^\alpha=
%\pa_t^{\alpha_0}
\pa_1^{\alpha_1}\cdots \pa_n^{\alpha_n}$
for a multi-index $\alpha=(%\alpha_0, 
\alpha_1, \ldots, \alpha_n)$.
% Similarly we write $\pa^\beta=\pa_t^{\beta_0}\pa_1^{\beta_1} \cdots \pa_n^{\beta_n}$
% with a multi-index $\beta=(\beta_0, \beta_1, \ldots, \beta_n)$.
We set
$$
\Gamma=\left(\Gamma_0, \Gamma_1, \ldots, \Gamma_{N}\right)
=\left(\pa_t, \pa_1, \ldots, \pa_n, (O_{ij})_{1\le i<j\le n}\right)
$$
with $N=n(n+1)/2$,
and $\Gamma^\beta=\Gamma_0^{\beta_0}\Gamma_1^{\beta_1}\cdots \Gamma_{N}
^{\beta_{N}}$ for a multi-index $\beta=(\beta_0,\beta_1, \ldots, \beta_N)$,
where $O_{ij}$ for $1\le i, j\le n$ is given by
$O_{ij}=x_i \pa_j-x_j \pa_i$.

For $r>0$ and $y\in \R^n$, $B_r(y)$ stands for an open ball of radius $r$ centered at $y$.
We write $B_r$ for $B_r(0)$.
Besides, we set $\Omega_r=\Omega \cap B_r$.

In what follows, for a constant $C$, when we write $C=C(p_1, \ldots, p_m)$
with $p_1, \ldots, p_m$ being some given constants or functions, 
it means that, with the space dimension $n$ and the obstacle ${\mathcal O}$ being fixed,
$C$ is a constant depending only on $p_1, \ldots, p_m$
(thus $C$ may depend also on $n$ and $\mathcal O$ actually). 

Then our main result reads as follows. 

\begin{theorem}\label{inverse}\
Let the space dimension $n$ be odd, and $n\ge 3$. Assume that 
${\mathcal O}$ is non-trapping, and ${\mathcal O} \subset B_1$.
Let $a\,(>1)$ be a fixed number. %$a>3$.
Then for any $\vec{f} \in X^\infty(\Omega)$ with 
$\supp \vec{f} \subset \overline{\Omega_a}$, there exists uniquely
$\vec{f}_+ \in \Hi^\infty(\R^n)$ 
%such that 
satisfying \eqref{ki3}.
%the following property:
Moreover, there exists a positive constant $\mu=\mu(a)$ having the following property: For any nonnegative integer $k$,
there exists a positive constant $C=C(k,a)$ such that
\begin{align}\label{iw0}
& \left\|\exp(\mu\jb{\,\cdot\,})\left( U(t)\vec{f}-U_0(t)\vec{f}_+ \right) \right\|_{\Hi^{k}(\Omega)}\\
& \qquad\qquad\qquad \le C\exp \left(-\mu t \right)\,\|\vec{f}\|_{\Hi^{k}(\Omega)} 
\quad \mbox{for} \ t \ge 0,
\nonumber\\
%\end{align}
%and 
%\begin{align} 
\label{iw01}
& \left\| \exp \left( %\frac{\sigma}{4} 
2\mu \jb{\, \cdot\, }\right)\,\vec{f}_+ \right\|_{{\mathcal W}^{\,k,\infty}(\R^n)} 
   \le C \|\vec{f}\|_{\Hi^{k+[n/2]+1}(\Omega)},
\end{align}
where $\jb{x}=\sqrt{1+|x|^2}$ for $x\in \R^n$,
and $[n/2]$ denotes the largest integer not exceeding $n/2$.
\end{theorem}

Theorem~\ref{inverse} will be proved in Section \ref{S3}.
Our proof of Theorem~\ref{inverse} relies on the exponential decay of
the local energy (see Lemma~\ref{LED} below), and this is the reason why
$n(\ge 3)$ is assumed to be odd and the obstacle ${\mathcal O}$ to be non-trapping.
For the notion of the non-trapping obstacle, we refer to
Melrose~\cite{Mel79} for instance (see also Shibata--Tsutsumi \cite{ShiTsu83, ShiTsu86}).
For example, star-shaped obstacles are known to be non-trapping.

Note that \eqref{iw01} implies that each component of
$\vec{f}_+$ belongs to the Schwartz class ${\mathcal S}$, 
the class of rapidly decreasing functions.
%%%%%%%%%% Original %%%%%%%%%%%%%%%%%%%
% Note that \eqref{iw01} implies that each component of
% $\vec{f}_+$ belongs to the Schwartz class ${\mathcal S}$, 
% the class of rapidly decreasing functions.
% As a by-product, the leading term of the perturbed wave is 
% deduced from (\ref{iw0}). 
% In fact, the leading term of the unperturbed wave with compactly
% supported data is well known as the Friedlander radiation field
% defined by
% \begin{align}\label{Radiation}
% {\mathcal F}_0[\vec{g}](s,\eta) =
% \frac{1}{2(2\pi)^{\frac{n-1}{2}}} 
% \sum_{j=0}^1 (-\pa_s)^{\frac{n-1}2-j} {\mathcal R}[g_j](s,\eta)
% \end{align}
% for $\vec{g}=(g_0, g_1)$.
% Here ${\mathcal R}[\varphi]$ denotes the Radon transform of $\varphi=\varphi(x)$, that is
% $$
% {\mathcal R}[\varphi](s,\eta)=
% \int_{y\cdot\eta=s} \varphi(y) dS_y,
% $$
% where $dS_y$ denotes the area element on the hyperplane $\{y;\, y\cdot\eta=s\}$.
% We remark that the same is true 
% for rapidly decreasing data (see Proposition~\ref{radiation0} below for the detail). 
% Since $\vec{f}_+$ enjoys a far better decay 
% estimate than we can expect in general for rapidly decreasing functions, 
% the rate of the convergence to the radiation field is given
% as follows:
%%%%%%%%%%%%%%%%%%%%%%%%%%%%%%%%%%%%%%%

Now we turn our attention to the asymptotic pointwise behavior of
the perturbed wave $U(t)\vec{f}$.
To describe the result, we define the Friedlander radiation field
${\mathcal F}_0[\vec{g}]$ by
\begin{align}\label{Radiation}
{\mathcal F}_0[\vec{g}](s,\eta) =
\frac{1}{2(2\pi)^{\frac{n-1}{2}}} 
\sum_{j=0}^1 (-\pa_s)^{\frac{n-1}2-j} {\mathcal R}[g_j](s,\eta)
\end{align}
for $\vec{g}=(g_0, g_1) \in \left(\mathcal{S}(\R^n)\right)^2$.
Here ${\mathcal R}[\varphi]$ denotes the Radon transform of $\varphi=\varphi(x)$, that is
$$
{\mathcal R}[\varphi](s,\eta)=
%2^{-1} (2\pi)^{-\frac{n-1}2} 
\int_{y\cdot\eta=s} \varphi(y) dS_y,
$$
where $dS_y$ denotes the area element on the hyperplane $\{y;\, y\cdot\eta=s\}$.
The radiation field ${\mathcal F}_0[\vec{g}]$ is introduced to describe the main part
of the unperturbed wave $U(t)\vec{g}$ for $\vec{g}\in \left(C^\infty_0(\R^3)\right)^2$
in Friedlander \cite{Fri62}. 
Lax--Phillips \cite{LaxPhi89} showed that the main part of the perturbed wave
can also be written in terms of the Friedlander radiation field of some function,
but the convergence rate seems not to have been obtained. Thus we would like to investigate
the convergence rate of the perturbed waves to the pointwise asymptotics described by
the Friedlander radiation field.
Our result is the following.
\begin{theorem}\label{radiation}\
Let the assumptions of Theorem \ref{inverse} hold.
Then for any $\vec{f} \in X^\infty(\Omega)$ with 
$\supp \vec{f} \subset \overline{\Omega_a}$,
there exists $\vec{f}_+ \in \left({\mathcal S}(\R^n)\right)^2$ satisfying
the following property:
For any nonnegative integer $k$, % and any $\nu\ge 0$, 
there exists a positive constant $C=C(k, a, \vec{f})$
such that, writing $x=r\omega$, for  $r\ge t/2\ge 1$ and 
$\omega=(\omega_1, \omega_2, \ldots, \omega_n)\in S^{n-1}$ we have
\begin{align}\label{ki6}
& \sum_{|\alpha|\le k} 
\bigl|
  \Gamma^\alpha  \bigl\{ 
                   u(t,x)-r^{-\frac{n-1}2} {\mathcal F}_0[\vec{f}_+](r-t,\omega)
                  \bigr\} 
\bigr| \\
& \qquad\qquad \le C (1+r+t)^{-\frac{n+1}2} %(1+|r-t|)^{-\nu}, 
\exp\left(-\frac{\mu}{2}|r-t|\right),
\nonumber\\
% \end{align}
% \begin{align}
\label{ki7}
& \sum_{|\alpha|\le k} \bigl|
                \Gamma^\alpha \bigl\{
                     \pa_t u(t,x)-(-1) r^{-\frac{n-1}2} 
                     \bigl(\pa_s {\mathcal F}_0[\vec{f}_+]\bigr)(r-t,\omega)
                               \bigr\} 
                              \bigr| \\
& \quad{}+\sum_{|\alpha|\le k}\sum_{j=1}^n \bigl|
                \Gamma^\alpha \bigl\{
                     \pa_j u(t,x)-\omega_j r^{-\frac{n-1}2} 
                     \bigl(\pa_s {\mathcal F}_0[\vec{f}_+]\bigr)(r-t,\omega)
                               \bigr\} 
                              \bigr| \nonumber \\
& \qquad\qquad \le C (1+r+t)^{-\frac{n+1}2} %(1+|r-t|)^{-\nu},
\exp\left( - \frac{\mu}{2} |r-t|\right),
\nonumber
\end{align}
where $u(t, x)$ is 
%the first component of $U(t)\vec{f}$, namely 
the solution to
\eqref{eq}--\eqref{id}
and $\mu=\mu(a)$ is the positive constant from Theorem \ref{inverse}.
%and the Friedlander radiation field ${\mathcal F}_0[\vec{g}]$ is given 
\end{theorem}

%%%%%%%%%%%%%%%%%%%%%%%%%%%%%%%%%%%%%%%%%%%%%%%%
% We note that the asymptotic behavior of $\pa_t u$ along the ray $x=(t+s)\omega$ 
% with $s$ being fixed was considered in Lax--Phillips \cite{LaxPhi89}, but
% the convergence rate seems not to have been obtained.
%%%%%%%%%%%%%%%%%%%%%%%%%%%%%%%%%%%%%%%%%%%%%%%%

%%%%%%%%%%%%%%%%%%%%%%%%%%%%%%%%%%%%%%%%%%%%%
The proof of Theorem \ref{radiation} will be given in Section \ref{ProofRadiation},
after obtaining the detailed convergence rate for the Cauchy problem in Section 4
(see Proposition \ref{radiation0}).
%%%%%%%%%%%%%%%%%%%%%%%%%%%%%%%%%%%%%%%%%%%%%

%\noindent
%\bf Remark.}
We underline that the decaying factor $\exp(-\mu|r-t|/2)$
in the above estimates is quite meaningful even if the 
initial data is compactly supported, say $\supp \vec{f} \subset \overline{\Omega_a}$,
unlike the case of the Cauchy problem.

In fact, the solution $u(t,x)$ for the mixed problem is 
identically zero for $r-t\ge a$ and $t \ge 0$,
in view of the domain of dependence (see Lemma~\ref{DoD} below).
On the other hand, it is 
not expected to vanish 
for $r-t\le -a$ in general, because of the presence of the obstacle.
Accordingly, the radiation field ${\mathcal F}_0[\vec{f}_+](s,\omega)$ 
for the solution to the mixed problem vanishes for $s \ge a$ and $\omega\in S^{n-1}$
due to \eqref{ki6}, although it is not supposed to be zero for 
$s \le -a$ and $\omega\in S^{n-1}$ in general.
In contrast to this, if there is no obstacle, it is known that 
the radiation field for compactly supported data vanishes also for $s\le -a$
(this property is closely connected to the Huygens principle; see Lemma~\ref{Huygens} below).

In conclusion, it is essential to extract the factor $\exp(-\mu |r-t|/2)$,
in order to describe the behavior for the mixed problem in the region $r-t \le -a$.
\medskip

% {\bf Remark.}
% Let $\supp \vec{f} \subset \overline{\Omega_a}$.
% Then, in view of the domain of dependence (see
% Lemma~\ref{DoD} below), we have $\supp (U(t)\vec{f}) \subset \overline{\Omega_{t+a}}$ for $t\ge 0$. Hence we have
% $$
% \lim_{t\to \infty} (t+s)^{\frac{n-1}{2}} u\left(t, (t+s)\omega\right)=0
% $$
% for any $s\ge a$ and any $\omega\in S^{n-1}$,
% where $u(t,\cdot)$ is the first component of $U(t)\vec{f}$.
% On the other hand, \eqref{ki6} yields
% $$
% \lim_{t\to\infty} (t+s)^{\frac{n-1}{2}} u\left(t, (t+s)\omega \right)={\mathcal F}_0[\vec{f}_+](s, \omega)
% $$
% for any $s\in \R$ and $\omega\in S^{n-1}$.
% Thus we conclude that
% $$
% {\mathcal F}_0[\vec{f}_+](s,\omega)=0
% $$
% for $s\ge a$ and $\omega\in S^{n-1}$. 
% However, because of the existence of the obstacle,
% ${\mathcal F}_0[\vec{f}_+](s, \omega)$ is not expected to vanish 
% for $s\le -a$ in general. This is the point where the decay factor $1+|r-t|$ plays its essential role. In contrast to this,
% if there is no obstacle, it is known that 
% the radiation field for compactly supported data vanishes also for $s\le -a$
% (this property is closely connected to the Huygens principle; see Lemma~\ref{Huygens} below).
% \medskip

%%%%%%%%%%%%%%%%%%%%%%%%%%%%%%%%%%
\section{Preliminaries}
%%%%%%%%%%%%%%%%%%%%%%%%%%%%%%%%%%%
%%%%%%%%%%%%%%%%%%%%%%%%%%%%%%
Let $Y$ be an open subset of $\R^n$, and $\Omega$ be as in the previous section.
For the notational convenience, we put 
\begin{align}\label{support}
& {\mathcal H}^\infty_a(Y)=\{\vec{f}=(f_0, f_1) \in {\mathcal H}^\infty(Y); 
\ \supp \vec{f} \subset \overline{Y\cap B_a} \},
\\
& X^\infty_a(\Omega)=\{\vec{f}=(f_0, f_1) \in X^\infty(\Omega); 
\ \supp \vec{f} \subset \overline{\Omega_a} \}
\end{align}
for $a>0$.

The following property is well known.
%%%%%%%%%%%%%%%%%%%%%%%%%%%%%%%%%%%%%%%%%%%%%%%%%
\begin{lemma}[Domain of dependence] 
\label{DoD}
Let $n$ be a positive integer.
Let $\tau, t_0 \in \R$ with $\tau<t_0$, and $x_0\in \R^n$.
We define 
$$
\Lambda(t_0, x_0, \tau)=\left\{(t,x)\in (\tau, t_0)\times \R^n;\, |x-x_0|< t_0-t \right\}.
$$ 
%and we write $B_r(x_0)$ for an open ball of radius $r$ centered at $x_0$.
Suppose that $\psi=\psi(t,x)$ satisfies
$$
(\pa_t^2-\Delta)\psi(t,x)=0, \quad (t,x)\in \Lambda(t_0, x_0, \tau).
$$
Then we have
\begin{equation}
\label{ConeEne}
\|\pa \psi(t)\|_{L^2\left( B_{t_0-t}(x_0) \right)}\le \|\pa \psi(\tau)\|_{L^2\left(B_{t_0-\tau}(x_0)\right)}, \quad t\in (\tau, t_0),
\end{equation}
where $\pa \psi=(\pa_t\psi, \nabla \psi)$. As a consequence, if we also assume
$$
\psi(\tau, x)=(\pa_t \psi)(\tau, x)=0, \quad x\in B_{t_0-\tau}(x_0),
$$
then we have $\psi(t,x)=0$ for any $(t,x)\in \Lambda(t_0, x_0, \tau)$.

%%%%%%%%%%%%%%%%%%%%%%%%%%%%%%%%%%%%%%%%%%%%%%%%%
The above assertions are also valid if we replace $\Lambda(t_0, x_0, \tau)$ by
$$
\Lambda^*(t_0, x_0, \tau)=\left\{(t,x)\in (2\tau-t_0, \tau)\times \R^n;\, |x-x_0|< t+t_0-2\tau \right\},
$$
and \eqref{ConeEne} by
$$
\|\pa \psi(t)\|_{L^2\left( B_{t+t_0-2\tau}(x_0) \right)}\le \|\pa \psi(\tau)\|_{L^2\left(B_{t_0-\tau}(x_0)\right)}, \quad t\in (2\tau-t_0, \tau).
$$
\end{lemma}
% Putting $\psi^*(t,x)=\psi(2\tau-t,x)$, we have $\pa_t^k\psi^*(t,x)=(-1)^k(\pa_t\psi)(2\tau-t,x)$.
% Hence 
From the lemma above, we see that  
$\vec{f}\in X_a^\infty(\Omega)$ (resp.~$\vec{g}\in \Hi^\infty_a(\R^n)$)
implies $\supp (U(t)\vec{f}) \subset \overline{\Omega_{|t|+a}}$
(resp.~$\supp \left(U_0(t)\vec{g}\right) \subset \overline{B_{|t|+a}}$).

In odd space dimensions, we have a stronger result.
\begin{lemma}[The Huygens principle]
\label{Huygens}
Let $n$ be an odd integer with $n\ge 3$.
Then $\vec{g}\in \Hi^\infty_a(\R^n)$ implies
$$
\supp \left(U_0(t) \vec{g}\right) \subset 
\left\{ x\in \R^n; |t|-a\le |x| \le |t|+a\right\},\quad t\in \R.
$$
\end{lemma}
This result follows immediately from the explicit expression of $U_0(t)\vec{g}$
(see \eqref{CouHil1} below).
%%%%%%%%%%%%%%%%%%%%%%%%%%%%%%%%%%%%%%%%%%%%%%%%%

Next we introduce 
the local energy decay of 
the perturbed wave at exponential rate (for the proof, see for instance 
Melrose \cite{Mel79}; see also Shibata--Tsutsumi \cite{ShiTsu83}). 

\begin{lemma}\label{LED}\
Let $n$ be odd and $n\ge 3$.
Assume that ${\mathcal O}$ is non-trapping, and ${\mathcal O} \subset B_1$.
Suppose that $a$, $b>1$, and $k$ is a nonnegative integer.
Then 
there exist two positive constants $C=C(k,a,b)$ and $\sigma=\sigma(a,b)$  such that  
for any $\vec{f} \in X^\infty_a(\Omega)$ we have
\begin{align}\label{ki5}
 \| U(t)\vec{f} \|_{\Hi^{k}(\Omega_b)} \le
C\exp(-{\sigma} t)\,\|\vec{f}\|_{\Hi^{k}(\Omega)} 
\quad \mbox{for} \ t \ge 0.
\end{align}
\end{lemma}

The following lemma, motivated by the arguments in Ikawa \cite{Ika00},
tells us that the perturbed wave can be decomposed
into the unperturbed wave and the correction term.
The former is the main part of the perturbed wave, while the latter 
takes care of the effect from the boundary and its size can be small compared
with the initial energy.
This lemma is crucial for proving Theorem \ref{inverse}.

\begin{lemma}\label{decomposition1}\
Let $n$, ${\mathcal O}$, and $a$ be as in 
Theorem \ref{inverse}.
Then, 
for any $\vec{f} \in X^\infty_a(\Omega)$ and $T(\ge a+2)$,
there exist $\vec{g}_1 \in \Hi_{T+a}^\infty(\R^n)$ and 
$\vec{f}_1 \in X^{\infty}_3(\Omega)$ 
satisfying
\begin{align}
 \label{iw1}
& U(t)\vec{f}=U_0(t-T)\vec{g}_1+U(t-T)\vec{f}_1,
\quad %\mbox{for} \ 
t \ge T, \\
   \label{iw3}
& \|\vec{g}_1\|_{\Hi^{k}(\R^n)} \le C_0 (1+T) \|\vec{f}\|_{\Hi^{k}(\Omega)},
\\ \label{iw21}
& \|\vec{f}_1\|_{\Hi^{k}(\Omega)} \le C_0 
     \exp(-\sigma T)\,\|\vec{f}\|_{\Hi^{k}(\Omega)}
\end{align}
for any nonnegative integer $k$
with some positive constants $C_0=C_0(k,a)$ and $\sigma=\sigma(a)$.
\end{lemma}

\begin{proof}
In this proof, various positive constants depending only on $k$ will be
indicated by the same $C_k$.
 
We put $T_0=T-2\,(\ge a)$.
If we set $\vec{\phi}=U(T_0)\vec{f}$, then 
$\vec{\phi} \in X^{\infty}(\Omega)$ and 
\begin{equation}
\label{A1}
\|\vec{\phi}\|_{\Hi^{k}(\Omega)} \le C_k (1+T_0) \|\vec{f}\|_{\Hi^{k}(\Omega)}
\end{equation}
for any nonnegative integer $k$.
Indeed, \eqref{A1} follows from the fact that we have
\begin{align}\label{iw31}
 \| U(t)\vec{f} \|_{\Hi^{k}(\Omega)} \le C_k (1+|t|) \|\vec{f}\|_{\Hi^{k}(\Omega)},
 \quad %\mbox{for} \ t \in \R.
 t\in \R
\end{align}
for any $\vec{f} \in X^\infty(\Omega)$. 
This estimate is a simple consequence of the energy estimate and
an elementary inequality
\begin{equation}\label{EleIneq}
\|v(t)\|_{L^2(\Omega)} \le \|v(t_0)\|_{L^2(\Omega)}+\int_{t_0}^t
\|\pa_t v(\tau)\|_{L^2(\Omega)} d\tau, \quad t\ge t_0,
\end{equation}
which is valid for any smooth function $v$.
%Here we denote the first component of $U(t)\vec{f}$ by $u(t)$.
Besides, in view of the domain of dependence (see Lemma~\ref{DoD}), we have $\supp \vec{\phi} \subset \overline{\Omega_{T_0+a}}$,
since $\supp \vec{f} \subset \overline{\Omega_a}$.

Next we extend $\vec{\phi}$ to $\vec{\psi} \in 
\Hi^{\infty}_{T_0+a}(\R^n)$ 
in such a way that 
$\vec{\psi}=\vec{\phi}$ in $\Omega$ and 
\begin{align} \label{iw4}
\|\vec{\psi}\|_{\Hi^{k}(\R^n)}
 \le C_k (1+T_0) \|\vec{f}\|_{\Hi^{k}(\Omega)}. 
\end{align}
To do this, we set $\vec{\phi}_0=\chi \vec{\phi}$ and 
$\vec{\phi}_\infty=(1-\chi) \vec{\phi}$, 
where $\chi$ is a smooth function on $\R^n$
satisfying $\chi(x)=1$ for $|x| \le 5$ and 
$\chi(x)=0$ for $|x| \ge 6$. 
Then $\vec{\phi}_0$ can be regarded as a function on a bounded domain 
$\Omega_6$, and we see from the Stein extension theorem
that $\vec{\phi}_0$ can be extended to
$\vec{\psi}_0 \in \Hi^{\infty}(\R^n)$ such that 
$\vec{\psi}_0=\vec{\phi}_0$ in $\Omega$ and 
\begin{equation}
\label{Stein01}
\|\vec{\psi}_0\|_{\Hi^{k}(\R^n)} \le C_k \|\vec{\phi}_0\|_{\Hi^{k}(\Omega)}
\end{equation}
for any nonnegative integer $k$ (refer to \cite{Ste70}). 
Recalling \eqref{A1} and setting
$\vec{\psi} := \vec{\psi}_0+\vec{\phi}_{\infty}$, we see that $\vec{\psi}$ has the desired properties.

Next we let $v$ be the solution of 
\begin{align}\label{eqC}
& (\partial_t^2-\Delta) v(t,x) =0, 
& (t,x) \in (T_0,\infty)\times \R^n,
\\ \label{idC}
& (v(T_0,x), (\partial_t v)(T_0,x))=\vec{\psi}(x),  
& x\in \R^n,
\end{align}
and we define $w=u-v$ in $[T_0,\infty)\times\Omega$, so that
\begin{align}\label{eqD}
& (\partial_t^2-\Delta) w(t,x) =0, & (t,x) \in (T_0,\infty) \times \Omega,
\\ \label{idD}
& w(T_0,x)=(\partial_t w)(T_0,x)=0, & x\in \Omega,
\end{align}
where $u(t,\cdot)$ denotes the first component of $U(t)\vec{f}$.
Furthermore, we define
\begin{align}\label{iw8}
& \vec{g}_1(x):=(v(T,x),(\partial_t v)(T,x))=U_0(2)\vec{\psi}(x),
\\ \label{iw9}
& \vec{f}_1(x):=(w(T,x),(\partial_t w)(T,x))=U(T)\vec{f}(x)-U_0(2)\vec{\psi}(x)
\end{align}
(recall $T=T_0+2$).
Then we easily get (\ref{iw3}) from (\ref{iw4}), because
we have
\begin{align}\label{iw91}
\| U_0(t) \vec{\psi}\|_{\Hi^{k}(\R^n)} \le 
C_k (1+|t|) \|\vec{\psi}\|_{\Hi^{k}(\R^n)},
\quad t\in \R
\end{align}
for $\vec{\psi} \in \Hi^\infty(\R^n)$.
This estimate is shown similarly to (\ref{iw31}).
Taking the domain of dependence into account, we have $\vec{g}_1\in{\mathcal H}^\infty_{T+a}(\R^n)$.

Next we consider $\vec{f}_1$.
Note that \eqref{eqD} and \eqref{idD} imply
\begin{align}\label{iw7}
w(t,x)=0 \quad \mbox{for} \ |x| \ge t-T_0+1,\ t \ge T_0
\end{align}
in view of the domain of dependence, because we have ${\mathcal O}\subset B_1$.
Hence $\supp \vec{f}_1 \subset \overline{\Omega_3}$, so that
\begin{align*}
\|\vec{f}_1\|_{\Hi^{k}(\Omega)} & \le 
\|U(T)\vec{f}\|_{\Hi^{k}(\Omega_3)}+\|U_0(2)\vec{\psi}\|_{\Hi^{k}(B_3)}
\\
& \le C \exp(-\sigma T)\,\|\vec{f}\|_{\Hi^{k}(\Omega)}+C_k
   \|\vec{\psi}\|_{\Hi^{k}(B_5)},
\end{align*}
thanks to (\ref{ki5}), \eqref{ConeEne} and \eqref{EleIneq}, where $C=C(k,a)$
and $\sigma=\sigma(a)$ are
positive constants.
Since $\vec{\psi}=\vec{\psi}_0$ in $B_5$, 
\eqref{Stein01} yields 
$$
 \|\vec{\psi}\|_{\Hi^{k}(B_5)} \le \|\vec{\psi}_0\|_{\Hi^{k}(\R^n)}
\le C_k \|\vec{\phi}_0\|_{\Hi^{k}(\Omega)}  \le C_k \|\vec{\phi}\|_{\Hi^{k}(\Omega_6)}.
$$ 
Recalling $\vec{\phi}=U(T_0)\vec{f}$ and using (\ref{ki5}) again, we obtain (\ref{iw21}).

In order to show that $\vec{f}_1 \in X^\infty_3(\Omega)$, 
it suffices to prove
\begin{align} \label{iw6}
w(t,x)=0 \quad \mbox{for} \ (t,x) \in [T_0+2,\infty) \times \partial \Omega.
\end{align}
Indeed, we already know $\vec{f}_1\in {\mathcal H}_3^\infty(\Omega)$;
as for the compatibility condition, writing $\vec{f}_1=(f_{1,0}, f_{1,1})$ and $f_{1,j}=\Delta f_{1, j-2}$ for $j\ge 2$,
we find $f_{1,j}(x)=(\pa_t^j w)(T,x)$ for $j\ge 0$, and \eqref{iw6} immediately leads
to $f_{1,j}=0$ on $\pa \Omega$ for $j\ge 0$.
Since $w=u-v$ and $\pa \Omega\subset B_1$, \eqref{iw6} is a consequence of (\ref{dc}) and
\begin{align}\label{iw5}
 v(t,x)=0 \quad \mbox{for} \ t \ge |x|+T_0+1.
\end{align}
To prove \eqref{iw5}, we define a function $z$ on $[0,\infty)\times \R^n$ by
\begin{align*}
z(t,x)=\begin{cases}
u(t,x) \  & \mbox{for} \ (t,x) \in [0,T_0] \times \Omega,
\\
0      \  & \mbox{for} \ (t,x) \in [0,T_0] \times \overline{\mathcal O},
\\
v(t,x) \  & \mbox{for} \ (t,x) \in (T_0,\infty) \times \R^n.
\end{cases}
\end{align*}
For $\ve>0$, let $\xi_\ve$ be a smooth function on $[0,\infty)\times \R^n$ such that
\begin{align*}
\xi_\ve(t,x)=\begin{cases}
1 \  & \mbox{for} \ |x| \ge 1 \ \mbox{or}\ t \ge T_0+\ve,
\\
0 \  & \mbox{for} \ (t,x) \in [0,T_0] \times {\mathcal O},
\end{cases}
\end{align*}
and $\xi_\ve(t,x)=\xi_\ve(0, x)$ for $(t,x) \in [0,T_0] \times \R^n$.
Then we have
\begin{align*}
& \supp\,\dal (\xi_\ve z)  \subset [0,T_0+\ve] \times \overline{B_1},
\\
& \supp\,(\xi_\ve z) (0,\cdot) \cup \supp\,\pa_t (\xi_\ve z) (0,\cdot) \subset \overline{B_{T_0}},
\end{align*}
since $a \le T_0$.
From the Duhamel principle, we have
$$
(\xi_\ve z)(t,\cdot)=U_0(t)\left((\xi_\ve z)(0), \pa_t(\xi_\ve z)(0)\right)
+\int_0^t U_0(t-\tau)\left(0, \dal(\xi_\ve z)(\tau)\right)d\tau.
$$
Thus by the Huygens principle (Lemma~\ref{Huygens})
we see that $v(t,x)=(\xi_\ve z)(t,x)=0$ for 
$t \ge |x|+T_0+\ve+1$, which implies (\ref{iw5}) because $\ve$ is arbitrary.

Finally, we prove (\ref{iw1}).
We see from (\ref{iw8}) and (\ref{iw9}) that (\ref{iw1}) holds at $t=T$.
Besides, for $(t,x) \in [T,\infty) \times \partial \Omega$ we have 
\begin{align*}
(U_0(t-T)\vec{g}_1)(x)+(U(t-T)\vec{f}_1)(x)
=&(U_0(t-T)\vec{g}_1 )(x) \\
=&(v(t,x), (\pa_t v)(t,x))=(0, 0)
\end{align*}
by (\ref{iw5}). 
It is apparent that we have $\dal(U_0(t-T)\vec{g}_1+U(t-T)\vec{f}_1)=0$ for $t\ge T$.
Hence we find (\ref{iw1}) by the uniqueness of the solution for the mixed problem. This completes the proof.
\end{proof}

\section{Proof of Theorem \ref{inverse}}
\label{S3}

In this section we prove Theorem \ref{inverse}.
The uniqueness is deduced from the following assertion:\
For given $\vec{f} \in X^\infty_a(\Omega)$, if $\vec{f}_+ \in \Hi^0(\R^n)$ 
satisfies
\begin{equation} \label{uni1}
\lim_{t \to \infty} \| U(t)\vec{f}-U_0(t)\vec{f}_+ \|_{\Hi_D(\Omega)}
=0,
\end{equation}
then $\vec{f}_+$ is determined uniquely.
To verify this assertion, suppose that $\vec{g}_+ \in \Hi^0(\R^n)$ also
satisfies $\displaystyle \lim_{t \to \infty} \| U(t)\vec{f}-U_0(t)\vec{g}_+ \|_{\Hi_D(\Omega)}
=0$, so that 
\begin{equation}
\label{unique01}
\lim_{t \to \infty} \| U_0(t)(\vec{f}_+ - \vec{g}_+) \|_{\Hi_D(\Omega)}=0.
\end{equation}
We also have 
\begin{equation}
\label{unique02}
\lim_{t \to \infty} \| U_0(t)(\vec{f}_+ - \vec{g}_+) \|_{\Hi_D(B_1)}=0.
\end{equation}
In fact, for any $\ve>0$, there exists $\vec{h} \in (C^\infty_0(\R^n))^2$ such that
$\| (\vec{f}_+ - \vec{g}_+) -\vec{h}\|_{\Hi_D(\R^n)} <\ve$.
Let $\supp \vec{h} \subset \overline{B_M}$.
Since the Huygens principle implies $U_0(t)\vec{h}=0$ for $|x| \le 1$ and $t \ge M+1$, we obtain
\begin{align*}
\| U_0(t)(\vec{f}_+ - \vec{g}_+) \|_{\Hi_D(B_1)}=
& \| U_0(t)(\vec{f}_+ - \vec{g}_+-\vec{h}) \|_{\Hi_D(B_1)}
\\
\le & \| \vec{f}_+ - \vec{g}_+-\vec{h} \|_{\Hi_D(\R^n)}
<\ve
\end{align*}
for $t \ge M+1$, which leads to \eqref{unique02}.
Here we have used the unitarity of $U_0(t)$ on $\Hi_D(\R^n)$.
From \eqref{unique01} and \eqref{unique02},
we see that 
$$
\|\vec{f}_+-\vec{g}_+\|_{\Hi_D(\R^n)}=\|U_0(t)(\vec{f}_+-\vec{g}_+)\|_{\Hi_D(\R^n)}
\to0 \quad (t\to \infty),
$$
which implies
$\vec{f}_+=\vec{g}_+$ in $\Hi_D(\R^n)$.
Since the H\"older inequality and the Sobolev imbedding theorem imply
that, for any $R>0$, there exists a positive constant $C_R$ such that we have
$$
 \|v\|_{L^2(B_R)} \le C_R \| v \|_{L^{2n/(n-2)}(\R^n)}
\le C_R \|\nabla v \|_{L^2(\R^n)}
$$
for any $v \in \dot{H}^1(\R^n)$, 
we conclude that $\vec{f}_+=\vec{g}_+$ in $\Hi^0(\R^n)$.

Next we consider the existence part.
We set $a_*=\max\{a, 3\}$, and we fix a nonnegative integer $k$.
We put $\mu=\sigma(a_*)/4$ and $C_1=C_0(k, a_*)$, where
$\sigma$ and $C_0$ are from Lemma \ref{decomposition1}.
We choose $T(\ge a_*+2)$ to be so large that 
$C_1 \exp(-\mu T)\le 1$.
Then we see from Lemma~\ref{decomposition1} that for $\vec{f} \in X^\infty_a(\Omega)$, 
there exist $\vec{g}_1 \in \Hi^{\infty}_{T+a_*}(\R^n)$ and
$\vec{f}_1 \in X^{\infty}_3(\Omega)$ satisfying (\ref{iw1}),
\begin{align*}
\|\vec{g}_1\|_{\Hi^{k}(\R^n)} \le & C_1(1+T) \|\vec{f}\|_{\Hi^{k}(\Omega)},\\
\intertext{and}
\|\vec{f}_1\|_{\Hi^{k}(\Omega)} \le & 
C_1 \exp(-4\mu T)\,\|\vec{f}\|_{\Hi^{k}(\Omega)}
\le \exp(-3\mu T)\,\|\vec{f}\|_{\Hi^{k}(\Omega)}.
\end{align*}
We apply Lemma~\ref{decomposition1} to $\vec{f}_1$ again to find 
$\vec{g}_2 \in \Hi_{T+a_*}^{\infty}(\R^n)$ and $\vec{f}_2 \in X^{\infty}_3(\Omega)$ 
for which we have
\begin{align*}
& U(t-T)\vec{f}_1=U_0(t-2T)\vec{g}_2+U(t-2T)\vec{f}_2
\quad \mbox{for} \ t \ge 2T,
\\ 
& \|\vec{g}_2\|_{\Hi^{k}(\R^n)} \le C_1 (1+T) \|\vec{f}_1\|_{\Hi^{k}(\Omega)}
  \le C_1(1+T)\exp(-3\mu T)\,\|\vec{f}\|_{\Hi^{k}(\Omega)},
\end{align*}
and
\begin{align*}
  \|\vec{f}_2\|_{\Hi^{k}(\Omega)} \le \exp(-3\mu T)\,\|\vec{f}_1\|_{\Hi^{k}(\Omega)}
  \le \exp(-6\mu T)\, \|\vec{f}\|_{\Hi^{k}(\Omega)}.
\end{align*}
Repeating the same procedure, we can construct sequences 
$\{\vec{g}_j\}_{j=1}^\infty \subset \Hi^{\infty}_{T+a_*}(\R^n)$ 
and $\{\vec{f}_j\}_{j=1}^\infty \subset X^{\infty}_3(\Omega)$ 
in such a way that 
\begin{align}\label{iw10}
 & U(t-(j-1)T)\vec{f}_{j-1} =U_0(t-jT)\vec{g}_j+U(t-jT)\vec{f}_j, \quad 
t \ge jT,
\\ \label{iw11}
  & \|\vec{g}_j\|_{\Hi^{k}(\R^n)} 
    \le C_1(1+T) \exp\left(-3\mu(j-1)T\right)\,\|\vec{f}\|_{\Hi^{k}(\Omega)},
\end{align}
and
\begin{align} \label{iw12}
&  \|\vec{f}_j\|_{\Hi^{k}(\Omega)} \le \exp(-3\mu j T) \,\|\vec{f}\|_{\Hi^{k}(\Omega)}
\end{align}
for $j\ge 1$, where we have put $\vec{f}_0=\vec{f}$.

Now, we define $\vec{f}_+=\sum_{j=1}^\infty U_0(-jT)\vec{g}_j$, which belongs to $\Hi^{\infty}(\R^n)$.
In fact, (\ref{iw91}) and (\ref{iw11}) lead to
\begin{align}\label{iw1a}
\|U_0(-jT)\vec{g}_j\|_{\Hi^{k}(\R^n)} & \le C (1+jT) 
  \exp\left(-3\mu(j-1) T\right)\,\|\vec{f}\|_{\Hi^{k}(\Omega)}
  \\ & \le C  \exp\left(-2\mu(j-1) T\right)\,\|\vec{f}\|_{\Hi^{k}(\Omega)},
  \nonumber
\end{align}
where $C$ is a constant depending on $k$ and $T$, but is independent of $j$.
Here we have used $(1+T+y)\exp(-\mu y)\le \mu^{-1} \exp\left(\mu(1+T)-1\right)$
for $y\in \R$.
Therefore we have
\begin{align}\label{iw13}
  \|\vec{f}_+\|_{\Hi^{k}(\R^n)} 
 \le \sum_{j=1}^\infty  C \left(\exp(-2\mu T)\right)^{j-1}\,\|\vec{f}\|_{\Hi^{k}(\Omega)}
 \le C \|\vec{f}\|_{\Hi^{k}(\Omega)}.
\end{align}

Next we prove (\ref{iw01}). 
For $\vec{h}=(h_0, h_1)$, we write 
\begin{equation}
\label{Not01}
|\vec{h}(x)|_{k}=\sum_{|\alpha|\le k+1}|\pa_x^\alpha h_0(x)|+\sum_{|\alpha|\le k} |\pa_x^\alpha h_1(x)|
\end{equation}
 in what follows.
Since $\supp \vec{g}_j\subset B_{T+a_*}$, 
the Huygens principle implies
\begin{equation}
 \label{support01}
 \supp \left(U_0(t-jT)\vec{g}_j \right)\subset \left\{x\in \R^n;\, 
 \bigl|\,|x|-|jT-t|\,\bigr|\le T+a_*\right\}
\end{equation}
for any natural number $j$ and $t\in \R$.
Hence it
follows from the Sobolev imbedding theorem and (\ref{iw1a}) that 
\begin{align}
\label{support02}
 |(U_0(-jT)\vec{g}_j)(x)|_{k} 
&  \leq C \|U_0(-jT)\vec{g}_j\|_{{\mathcal H}^{k+[n/2]+1}(\R^n)}
\\
&  \leq C \exp\left(-2\mu (j-1)T \right) \|\vec{f}\|_{{\mathcal H}^{k+[n/2]+1}(\Omega)}
\nonumber\\
&  \leq C \exp\left(-2\mu |x| \right)
                  \|\vec{f}\|_{{\mathcal H}^{k+[n/2]+1}(\Omega)}
\nonumber
\end{align}
for $x\in \supp \left(U_0(-jT)\vec{g}_j\right)$,
where $C$ is a constant depending on $k$, $a$ and $T$, but is independent of $j$ and $x$.
Noting that, for each fixed $x\in \Omega$, the number of $j$ 
for which we have $x\in \supp \left(U_0(-jT)\vec{g}_j\right)$ is at most $[2(T+a_*)/T]+1$ (cf.~\eqref{support01}),
we obtain (\ref{iw01}) from \eqref{support02}.

%%%%%%%%%%%%%%%%
Next we prove (\ref{iw0}). 
For $t \ge T$, we find a positive integer $J$
such that $t \in [JT,(J+1)T)$.
By (\ref{iw10}) with $j=1,\dots,J$ we have
$$
U(t)\vec{f}=\sum_{j=1}^J U_0(t-jT)\vec{g}_j+U(t-JT)\vec{f}_J.
$$
Since $U_0(t)\vec{f}_+=\sum_{j=1}^\infty U_0(t-jT)\vec{g}_j$, we get
\begin{align}
\label{A2}
&  \| e^{\mu\jb{\,\cdot\,}} (U(t)\vec{f}-U_0(t)\vec{f}_+) \|_{\Hi^{k}(\Omega)} 
\\
& \ \le \sum_{j=J+1}^\infty \|e^{\mu\jb{\,\cdot\,}}U_0(t-jT)\vec{g}_j\|_{\Hi^{k}(\Omega)} 
+\|e^{\mu\jb{\,\cdot\,}} U(t-JT)\vec{f}_J\|_{\Hi^{k}(\Omega)}.
\nonumber
\end{align}
Note that \eqref{A2} is also valid for $0\le t<T$, by regarding $J=0$ and
$\vec{f}_0=\vec{f}$. So we assume $J\ge 0$ and $t\in [JT, (J+1)T)$ 
in the following.
Since $|t-JT| \le T$ for $t \in [JT,(J+1)T)$, 
we get
$$
\sum_{|\alpha|\le k+1} \left|\pa_x^\alpha \exp(\mu \jb{x})\right| \le C_k \exp\left(\mu \jb{2T+a_*}\right)
$$
for $x\in \supp (U(t-JT)\vec{f}_J)$ with some positive constant $C_k$ depending only on $k$
and $\mu(=\sigma(a_*)/4)$.
Thus the second term on the right-hand side of \eqref{A2}
is estimated 
by 
\begin{align*}
 C (1+|t-JT|) \|\vec{f}_J\|_{\Hi^{k}(\Omega)}
\le C (1+T) \exp(-3\mu JT)\,\|\vec{f}\|_{\Hi^{k}(\Omega)}.
\end{align*}
Here we have used (\ref{iw31}) and (\ref{iw12}), % for the last inequality 
and the constant $C=C(k, a, T)$ is independent of $J$.
From \eqref{support01}, we get
\begin{align*}
\sum_{|\alpha|\le k+1}
\left|\pa_x^\alpha \exp (\mu \jb{x}) \right| \le & C_k \exp \left(\mu(1+T+a_*+jT-t)\right) 
\\ \le & C_k \exp\left(\mu(1+T+a_*+jT)\right)
\end{align*}
for $x\in \supp\,(U_0(t-jT)\vec{g}_j)$ with $j\ge J+1$, where $C_k$ is a positive constant
depending only on $k$ and $\mu(=\sigma(a_*)/4)$.
Hence, it follows from (\ref{iw91}) and (\ref{iw11}) that
\begin{align*}
 \|e^{\mu\jb{\,\cdot\,}}U_0(t-jT)\vec{g}_j\|_{\Hi^{k}(\Omega)} & \le Ce^{\mu j T}(1+|t-jT|)
    \|\vec{g}_j\|_{\Hi^{k}(\R^n)}
\\
& \le  C(1+jT)  e^{\mu jT-3\mu(j-1) T}\|\vec{f}\|_{\Hi^{k}(\Omega)}
\\
& \le  C e^{-\mu(j-1)T}\|\vec{f}\|_{\Hi^{k}(\Omega)}
\end{align*}
for $j \ge J+1$ and $t \in [JT,(J+1)T)$, where $C$ is a constant independent of $j$ and $J$.
Thus the first term on the right-hand side of \eqref{A2} is evaluated by
$C \exp(-\mu JT)\,\|\vec{f}\|_{\Hi^{k}(\Omega)}$, where $C$ is a constant independent of $J$.
Therefore, (\ref{iw0}) holds for $t \in [JT,(J+1)T)$ with $J\ge 0$,
and hence for all %$t \ge T$.
$t\ge 0$.

% On the other hand, for $t \in [0,T)$, we have from (\ref{iw31}) and (\ref{iw91})
% $$
% \|e^{\mu\jb{\,\cdot\,}} (U(t)\vec{f}-U_0(t)\vec{f}_+) \|_{\Hi^{k}(\Omega)} \le 
% C_k (1+t) \{\|\vec{f}\|_{\Hi^{k}(\Omega)}
% +\|\vec{f}_+\|_{\Hi^{k}(\R^n)}\}, 
% $$
% which gives (\ref{iw0}) for $t \in [0,T)$ by (\ref{iw13}).

Finally, we remark that by the uniqueness result, 
$\vec{f}_+$ being constructed in the above is independent of $k$, although the construction itself depends on $k$
through the choice of $T$.
%%%%%%%%%%%%
In fact, let $\vec{f}_+^{\, (1)}$ and $\vec{f}_+^{\, (2)}$ denote
$\vec{f}_+$ constructed in the above with the choice of $k=k_1$ and
$k=k_2$, respectively, where $k_1$ and $k_2$ are nonnegative integers.
%%%%%%%%%%%%
Then, from \eqref{iw0}, \eqref{uni1} is valid 
for $\vec{f}_+ = \vec{f}_+^{\, (1)}$ and
$\vec{f}_+=\vec{f}_+^{\, (2)}$. Hence the uniqueness of $\vec{f}_+$ satisfying
\eqref{uni1} implies $\vec{f}_+^{\, (1)}=\vec{f}_+^{\, (2)}$.

This completes the proof of Theorem \ref{inverse}.
\qed
%%%%%%%%%%%%%%%%%%%%%%%%%%%%%%%%%%%%%%%
%\newpage
%%%%%%%%%%%%%%%%%%%%%%%%%%%%%%%%%%%%%%%
\section{%Appendix: 
The Friedlander Radiation Field for Rapidly Decreasing Data}
%%%%%%%%%%%%%%%%%%%%%%%%%%%%%%%%%%%%%%%
Our aim in this section is to discuss the Friedlander radiation field
for the Cauchy problem with rapidly decreasing data.
The case of compactly supported data is well known (see Friedlander \cite{Fri62, Fri64, Fri67};
see also H\"ormander \cite{Hoe97} and John \cite{Joh}).
The case of rapidly decreasing data was also treated in \cite{Hoe97}
through the conformal compactification of the Minkowski space.
But the decay away from the light cone was neglected there. 
Hence we would like to obtain a more detailed estimate, restricting
our attention to the odd space dimensional case.

As is known, the behavior of the solution away from the cone is 
closely related to the decay property of the data. % at spatial infinity.
Because the scattering data $\vec{f}_+$ obtained in Theorem \ref{inverse}
satisfies the stronger decay property than general functions in ${\mathcal S}(\R^n)$,
we introduce the following class of the data.
%We would like to also treat rapidly decreasing data with some additional decay factor.
Throughout this section, $\chi=\chi(s)$ is some given non-decreasing function 
of $s\ge 0$, satisfying $\chi(s)\ge 1$ for all $s\ge 0$.
For $\varphi\in C^\infty(\R^n)$, $m\ge 0$ and a nonnegative integer $k$, we define
$$
\|\varphi\|_{\chi,k,m}=\left(\sup_{x\in \R^n} \sum_{|\alpha|\le k}  
(1+|x|^2)^m \chi^2(|x|)
\left|\pa_x^\alpha \varphi(x)\right|^2\right)^{1/2},
$$
and let ${\mathcal S}_\chi(\R^n)$ be the set of all $\varphi\in C^\infty(\R^n)$ 
satisfying $\|\varphi\|_{\chi, m,k}<\infty$ for any nonnegative integers $m$ and $k$.
Apparently we have $\mathcal{S}_\chi(\R^n)\subset {\mathcal S}(\R^n)$,
where ${\mathcal S}(\R^n)$ is the Schwartz class, the set of rapidly decreasing functions.
Note that ${\mathcal S}_\chi(\R^n)={\mathcal S}(\R^n)$ if $\chi$ is identically equal to $1$. Our main result in this section is the following.
\begin{proposition}\label{radiation0}\
Let $n$ be an odd integer with $n\ge 3$, and let $\nu \ge 0$.
For any $\vec{f} \in \left({\mathcal S}_\chi (\R^n)\right)^2$ and any
multi-index $\alpha$, there exists a positive constant $C=C(\alpha, \nu, \vec{f})$ 
such that we have
\begin{align}\label{ki60}
& \left| \Gamma^\alpha\left\{
u(t,x)-r^{-\frac{n-1}2} {\mathcal F}_0[\vec{f}](r-t, \omega) \right\} 
\right| \\
& \qquad\qquad \le
C (1+t+r)^{-\frac{n+1}2}(1+|r-t|)^{-\nu}\chi^{-1}(|r-t|), \nonumber\\
%\end{align}
%\begin{align}
\label{ki70}
& \left|
\Gamma^\alpha\left\{
\pa_t u(t,x)-(-1) r^{-\frac{n-1}2} (\pa_s{\mathcal F}_0[\vec{f}])(r-t,\omega)
\right\}
\right|\\
& \quad {}+\sum_{j=1}^n \left|
\Gamma^\alpha\left\{
\pa_j u(t,x)-\omega_j r^{-\frac{n-1}2} (\pa_s{\mathcal F}_0[\vec{f}])(r-t,\omega)
\right\}
\right| \nonumber\\
&\qquad\qquad \le C (1+t+r)^{-\frac{n+1}2}(1+|r-t|)^{-\nu}\chi^{-1}(|r-t|) \nonumber
\end{align}
for $r\ge t/2\ge 1$ with $r=|x|$ and $\omega=(\omega_1,\ldots, \omega_n)=r^{-1}x$, 
where $u(t, \cdot)$ is the first component of $U_0(t)\vec{f}$, and
the radiation field ${\mathcal F}_0[\vec f](s, \eta)$ is given by \eqref{Radiation}.
%as in Theorem~\ref{radiation}.
% $\pa^\alpha$ denotes $\pa_t^{\alpha_0}\pa_1^{\alpha_1}\cdots \pa_n^{\alpha_n}$
% for a multi-index $\alpha=(\alpha_0, \alpha_1, \ldots, \alpha_n)$.
\end{proposition}
%%%%%%%%%%%%%%%%%%%%%%%%%%%%%%%%%%%%%
% For the sake of completeness, 
We will give a proof of this proposition, taking a fundamental approach based on the explicit representation of $U_0(t)\vec{f}$,
instead of using the conformal compactification.

% For $\varphi\in {\mathcal S}(\R^n)$, and nonnegative integers $m$ and $k$, we define
% $$
% \|\varphi\|_{m,k}=\left(\sup_{x\in \R^n} \sum_{|\alpha|\le k}  (1+|x|^2)^m 
% \left|\pa_x^\alpha \varphi(x)\right|^2\right)^{1/2}.
% $$
First we state some basic properties of the Radon transform.
We recall that the Radon transform ${\mathcal R}[\varphi](s, \eta)$ for
$\varphi\in {\mathcal S}(\R^n)$ is
defined by
$$
{\mathcal R}[\varphi](s, \eta)=\int_{\Pi(s,\eta)} \varphi(y) dS_y, \quad (s,\eta)\in \R\times S^{n-1},
$$
where $\Pi(s,\eta)=\{y\in \R^n; y\cdot\eta=s\}$, and $dS_y$ denotes
the area element on $\Pi(s,\eta)$.
For $\eta\in S^{n-1}$ and a smooth function $\varphi=\varphi(y)$ on $\R^n$,
$D_\eta \varphi$ denotes the directional derivative of $\varphi$ in the direction $\eta$;
in other words, we define 
$(D_\eta \varphi)(y)=\eta \cdot \nabla_y \varphi(y)$.
We write
$$
o_{ij}=\eta_i\pa_{\eta_j}-\eta_j\pa_{\eta_i},
%\ O_{ij}=y_i\pa_{y_j}-y_j\pa_{y_i}, 
\quad 1\le i, j\le n.
$$
We put $o=(o_1,\ldots, o_{n(n-1)/2})=(o_{ij})_{1\le i<j\le n}$, 
where $o_{ij}$'s are regarded to be arranged in dictionary order.
We write $o^\alpha=o_1^{\alpha_1}\cdots o_{d}^{\alpha_d}$
with a multi-index $\alpha$, where $d=n(n-1)/2$.
$O^\alpha$ is similarly defined using $O_{ij}$ instead of $o_{ij}$, where
$(O_{ij}\varphi)(y)=y_i (\pa_j\varphi)(y)-y_j(\pa_i \varphi)(y)$ as before.

It is easy to check
\begin{align}
\label{RadonDeri01}
\pa_s {\mathcal R}[\varphi](s,\eta)=& {\mathcal R}[D_\eta \varphi](s,\eta)
\left(=\int_{\Pi(s,\eta)} (D_\eta \varphi)(y) dS_y
\right), \\
\label{RadonDeri02}
o_{ij} {\mathcal R}[\varphi](s,\eta)=& {\mathcal R}[O_{ij}\varphi](s, \eta),
\quad 1\le i<j\le n
\end{align}
for $\varphi\in {\mathcal S}(\R^n)$.
Because integrals over $\Pi(s,\eta)$ of
directional derivatives of $\varphi$ in directions proportional to $\Pi(s,\eta)$ vanish,
we get
\begin{align}
\label{RadonDeri03}
{\mathcal R}[\pa_i\varphi](s,\eta)=&
{\mathcal R}[\eta_i D_\eta \psi](s,\eta)= \eta_i\pa_s
{\mathcal R}[\varphi](s,\eta) 
\end{align}
for $1\le i\le n$.

%%%%%%%%%%%%%%%%%
We observe that if $\varphi\in {\mathcal S}_\chi(\R^n)$, then we have 
%${\mathcal R}[\varphi]\in {\mathcal S}(\R\times S^{n-1})$ in the sense that we have
%\begin{equation}\label{DecayRadon}
%\sup_{(s,\theta)\in \R\times S^{n-1}}
%(1+s^2)^{\frac{N}{2}}\chi(s) \left|\pa_s^j o^\alpha {\mathcal R}[\varphi](s,\theta)\right|
%<\infty 
%\end{equation}
\begin{equation}
\label{DecayRadon02}
\left|\pa_s^jo^\alpha {\mathcal R}[\varphi](s, \eta)\right|
\le C_{j,\alpha} \|\varphi\|_{\chi, j+|\alpha|, \mu+n+|\alpha|} (1+s^2)^{-\frac{\mu}{2}}\chi^{-1}(|s|)
\end{equation}
for any $(s,\eta)\in \R\times S^{n-1}$, any $\mu\ge 0$,
any nonnegative integer  $j$,
and for any multi-index $\alpha$.
Here $C_{j,\alpha}$ denotes a positive constant depending only on $j$ and $\alpha$.
In fact, writing $\rho=|y-(y\cdot\eta)\eta|$,
we have $|y|^2=s^2+\rho^2$ for $y\in \Pi(s,\eta)$.
Hence we get
\begin{align}
\label{DecayRad02a}
& \left|(D^j_\eta O^\alpha \varphi)(y)\right| \\
& \quad \le C_{j,\alpha}(1+s^2+\rho^2)^{-\frac{\mu}{2}}(1+\rho)^{-n}\chi^{-1}(|s|)
\|\varphi\|_{\chi, j+|\alpha|, \mu+n+|\alpha|}
\nonumber
\end{align}
for $y\in \Pi(s,\eta)$.
In view of \eqref{RadonDeri01} and \eqref{RadonDeri02}, we find \eqref{DecayRadon02}.

We also notice that if $\varphi\in {\mathcal S}_\chi(\R^n)$, then we have 
\begin{align}
\label{RadonDeri04}
& \biggl|
    \pa_x^\alpha O^\beta
         \left\{
          (\pa_s^k {\mathcal R}[\varphi])\left(|x|-t, \frac{x}{|x|}\right) 
         \right\}\\
& \qquad\qquad\qquad {}-(\pa_s^k {\mathcal R}[\pa_x^\alpha O^\beta \varphi])\left(|x|-t, \frac{x}{|x|}\right)
\biggr|
\nonumber\\ 
& \quad \le C\frac{\|\varphi\|_{\chi, k+|\alpha|+|\beta|, \nu+n+|\alpha|+|\beta|}}{(1+t+|x|)\left(1+\bigl| |x|-t\bigr|\right)^{\nu}\chi\left(\bigl| |x|-t \bigr| \right)}
% & \quad \le C\|\varphi\|_{k+|\alpha|, \nu+n+|\alpha|}(1+t+|x|)^{-1}\left(1+|\,|x|-t|\right)^{-\nu}
\nonumber
\end{align}
for $|x|\ge t/2\ge 1$, $\nu\ge 0$, any nonnegative integer $k$, and
any multi-indices $\alpha$, $\beta$, where $C=C(k,\nu, \alpha, \beta)$ is a positive constant.
Since for any $\psi \in C^\infty (S^{n-1})$ we have $O^\beta \{ \psi(|x|^{-1}x)\} = (o^\beta \psi) (|x|^{-1}x)$, 
it suffices to show \eqref{RadonDeri04} for $\beta=0$, thanks to \eqref{RadonDeri02}.
By \eqref{RadonDeri03} we have 
\begin{align*}
& \pa_i\left\{
(\pa_s^k {\mathcal R}[\varphi])(|x|-t, |x|^{-1} x) 
\right\} \\
& = \left. \left(\eta_i 
\pa_s^{k+1}{\mathcal R}[\varphi](s,\eta) -|x|^{-1} \sum_{j=1}^n \eta_j o_{ij}\pa_s^k{\mathcal R}[\varphi](s,\eta)\right)\right|_{(s,\eta)=(|x|-t, |x|^{-1}x)}\\
& = (\pa_s^k {\mathcal R}[\pa_i\varphi])(|x|-t, |x|^{-1}x)-\sum_{j=1}^n \frac{x_j}{|x|^2}
( o_{ij} \pa_s^k {\mathcal R}[\varphi])(|x|-t, |x|^{-1}x)
\end{align*}
for $1\le i\le n$ and any nonnegative integer $k$.
Therefore, \eqref{DecayRadon02} implies
\begin{align} \nonumber
& \left|\pa_i\left\{
(\pa_s^k {\mathcal R}[\varphi])(|x|-t, |x|^{-1}x) 
\right\}
 {}-(\pa_s^k {\mathcal R}[\pa_i \varphi])(|x|-t, |x|^{-1} x)\right|\\
& \quad \le C\frac{\|\varphi\|_{\chi, k+1, \nu+n+1}}{(1+t+|x|)
\left(1+\bigl||x|-t\bigr|\right)^{\nu}\chi\left(\bigl||x|-t\bigr|\right)}
\nonumber
% & \ \le C\|\varphi\|_{\chi, k+1, \nu+n+1}(1+t+|x|)^{-1}
% \left(1+|\,|x|-t|\right)^{-\nu}\chi^{-1}(|\,|x|-t|)
% \nonumber
\end{align}
for $|x|\ge t/2\ge 1$, $1\le i\le n$, and $\nu\ge 0$,
where $C$ is a positive constant depending on $k$ and $\nu$. 
Hence \eqref{RadonDeri04} holds for $|\alpha|=1$.
Similarly we obtain it for general $\alpha$.

%%%%%%%%%%%%%%%%%%%%%%%%%%%%%%%%%%%%%%%
We now turn our attention to the explicit representation of 
$u(t, \cdot)$, the first component of $U_0(t)\vec{f}$.
Let $n(\ge 3)$ be an odd integer.
It is known that when $\vec{f}=(0,\varphi)$, $u(t,x)$ is expressed by
the following integral\,:
\begin{equation}\label{CouHil1}
E[\varphi](t,x)=\frac{\sqrt{\pi}}{2\Gamma(n/2)}
 \left( \frac{1}{2t} \frac{\pa}{\pa t} \right)^{\frac{n-3}2}
 (t^{n-2} Q[\varphi](t,x)),
\end{equation}
where $\Gamma(s)$ is the Gamma function and we put
\begin{equation}
\label{CouHil2}
Q[\varphi](t,x)=\frac1{A_n} \int_{\theta\in S^{n-1}} 
 \varphi(x+t \theta) dS_{\theta}^\prime
\end{equation}
for $\varphi \in {\mathcal S}(\R^n)$ and $(t,x)\in (0,\infty)\times \R^n$.
Here $A_n$ is the total measure of $S^{n-1}$, that is $A_n=2 \pi^{n/2}/\Gamma(n/2)$,
and $dS'_\theta$ is the area element on $S^{n-1}$
(see, e.g., Courant and Hilbert \cite[Chapter VI, Section 12]{CouHil62}).
Therefore, in general, $u(t, x)$ can be written as
\begin{equation*}
u(t,x)=\pa_t E[f_0](t,x)+E[f_1](t,x).
\end{equation*} 
We also have
\begin{equation}
\label{FundamentalSol02}
\pa_t^k\pa_x^\alpha O^\beta u(t,x)=\pa_t^{k+1} E[\pa_x^\alpha O^\beta f_0](t,x)+\pa_t^k E[\pa_x^\alpha O^\beta f_1] (t,x)
\end{equation}
for any nonnegative integer $k$ and any multi-indices $\alpha$, $\beta$
(note that $\varphi\in {\mathcal S}_\chi(\R^n)$ implies $\pa_x^\alpha O^\beta\varphi\in {\mathcal S}_\chi(\R^n)$
for any multi-indices $\alpha$ and $\beta$).
Hence, once we establish that there exist a large integer $N$ and a positive constant
$C=C(k)$ such that 
\begin{align}\label{CouHil3}
& \left|\pa_t^k E[\varphi](t, x)-\frac{1}{2(2\pi r)^{\frac{n-1}2}} 
\left((-\pa_s)^{\frac{n-3}2+k} {\mathcal R}[\varphi]\right) (r-t, \omega)\right|
\\ \nonumber
& \quad \le C \|\varphi\|_{\chi, \frac{n-1}2+k, N}(1+t+r)^{-\frac{n+1}2}(1+|r-t|)^{-\nu}\chi^{-1}(|r-t|)
\end{align}
holds for $\varphi\in {\mathcal S}_\chi(\R^n)$, $\nu \ge 0$, $r(=|x|)\ge t/2\ge 1$, and
$\omega=r^{-1}x$, 
we can conclude that Proposition~\ref{radiation0} is valid, 
in view of \eqref{Radiation}, \eqref{RadonDeri03},
\eqref{DecayRadon02}, \eqref{RadonDeri04} and \eqref{FundamentalSol02}.

In order to prove \eqref{CouHil3}, we observe that
\begin{align}
\label{CouHil4}
\pa_t^k E[\varphi](t,x)= \frac{\sqrt{\pi}}{2^{\frac{n-1}2} \Gamma(n/2)}
 \sum_{\ell=0}^{\frac{n-3}2+k} a_\ell\,t^{1-k+\ell} \pa_t^\ell Q[\varphi](t,x)
\end{align}
for $k\ge 0$, where $a_l$ are suitable constants with $a_\ell=1$ for $\ell=\frac{n-3}2+k$.
Let $\varphi\in {\mathcal S}_\chi(\R^n)$, $\nu \ge 0$, $r(=|x|)\ge t/2\ge 1$, and
$\omega=r^{-1}x$ in the following.
Then we see that \eqref{CouHil3} follows from the estimate
\begin{align}\label{CouHil5}
& \left|\pa_t^\ell Q[\varphi](t, x)-\frac{1}{A_n t^{n-1}} 
\left((-\pa_s)^{\ell} {\mathcal R}[\varphi]\right) (r-t, \omega)\right|
\\ \nonumber
& \qquad \le C \|\varphi\|_{\chi, \frac{n-1}2+k, N}(1+t+r)^{-n}(1+|r-t|)^{-\nu}\chi^{-1}(|r-t|).
\end{align}
In fact, when $2t \ge r\ge t/2\ge 1$, since we have
\begin{align*}
 \left|t^{-\frac{n-1}2}-r^{-\frac{n-1}2}\right|
 \le C(1+t+r)^{-\frac{n+1}2}(1+|r-t|),
\end{align*}
from \eqref{DecayRadon02} we see that it suffices to prove \eqref{CouHil3} 
with $(2\pi r)^{\frac{n-1}{2}}$ in its left-hand side being replaced by $(2\pi t)^{\frac{n-1}{2}}$.
This replaced estimate can be easily proved by \eqref{CouHil4}, \eqref{CouHil5} and \eqref{DecayRadon02}.
On the other hand, when $r>2t$, from \eqref{DecayRadon02} we get
$$
 \left|\left((-\pa_s)^\ell {\mathcal R}[\varphi]\right)(r-t, \omega)\right|
 \le C \|\varphi\|_{\chi, \ell, \mu+n} (1+t+r)^{-\mu}\chi^{-1}(|r-t|)
$$
for any $\mu\ge 0$.
Hence, using  \eqref{CouHil4} and \eqref{CouHil5},  we find that both  $|\pa_t^kE[\varphi](t,x)|$ and $|(-\pa_s)^{\frac{n-3}{2}+k} {\mathcal R}[\varphi])(r-t,\omega)|$ are bounded from above by the right-hand side
of \eqref{CouHil3}. Now we have seen that our task is to prove \eqref{CouHil5}. 

It follows from \eqref{CouHil2} that
\begin{align}\label{CouHil6}
\pa_t^\ell Q[\varphi](t,x) & =\frac1{A_n} \int_{S^{n-1}} 
 (D_\theta^\ell \varphi)(x+t \theta) dS_{\theta}^\prime
\\ \nonumber
& =\frac1{A_n} \sum_{|\alpha|=\ell} \int_{S^{n-1}} 
 c_{\ell,\alpha}(\theta) (\pa^\alpha  \varphi)(x+t \theta) dS_{\theta}^\prime
\\  \nonumber
& =\frac1{A_n t^{n-1}} \sum_{|\alpha|=\ell} 
\int_{S(t,x)} c_{\ell,\alpha} \left(t^{-1}(y-x)\right)\pa_y^\alpha \varphi(y) dS_y^*
\end{align}
with some polynomial $c_{\ell, \alpha}$, where 
$S(t, x)=\{y\in \R^n; |y-x|=t\}$ and $dS_y^*$ stands for the area element on $S(t,x)$
(recall that $D_\theta \varphi$ is the directional derivative of $\varphi$ in the direction $\theta$).
On the other hand, \eqref{RadonDeri01} implies that
\begin{equation}
 (-\pa_s)^\ell {\mathcal R}[\varphi](r-t, \omega)=\sum_{|\alpha|=\ell} {c_{\ell,\alpha}(-\omega)\mathcal R}[\pa_y^\alpha \varphi](r-t,\omega).
\end{equation}
Thus our task of proving \eqref{CouHil5} is reduced to the estimate
\begin{align}\label{CouHil7}
& \biggl|
\int_{S(t,x)} c_{\ell,\alpha} \left(t^{-1}(y-x)\right)\pa_y^\alpha \varphi(y) dS_y^*\\
& \qquad\qquad\qquad\qquad\qquad
{}-c_{\ell,\alpha}(-\omega) {\mathcal R}[\pa_y^\alpha \varphi](r-t, \omega) \biggr|
\nonumber\\ 
\nonumber
& \qquad \le C \|\varphi\|_{\chi, \frac{n-1}2+k, N}(1+t+r)^{-1}(1+|r-t|)^{-\nu}\chi^{-1}(|r-t|).
\end{align}

In order to proceed further, we decompose the integral over $S(t,x)$ as follows.
Let $\varepsilon$ be a small and positive constant.
For $r>0$, $t>0$ and $\omega \in S^{n-1}$, we set
\begin{align*}
\Lambda_\varepsilon^1(t, r, \omega)=& \{y\in S(t,r \omega);  |y|>(t+r)^{\varepsilon} \}, \\
\Lambda_\varepsilon^2(t, r, \omega)=& \{y\in S(t,r \omega); |y|\le (t+r)^{\varepsilon} \}.
\end{align*}
When $|r-t|> (t+r)^\varepsilon$, we have $S(t, r \omega)=\Lambda_\varepsilon^1(t, r, \omega)$.
Therefore, using \eqref{DecayRadon02} and Lemma~\ref{radiation1} below
to estimate $\int_{S(t,x)} c_{\ell,\alpha} \left(t^{-1}(y-x)\right)\pa_y^\alpha \varphi(y) dS_y^*$
and $c_{\ell,\alpha}(-\omega){\mathcal R}[ \pa_y^\alpha \varphi](r-t, \omega)$,
respectively, we obtain \eqref{CouHil7}.
On the other hand, when $|r-t|\le (t+r)^\varepsilon$,
\eqref{CouHil7} is a consequence of Lemmas~\ref{radiation1} and \ref{radiation2} below.

%%%%%%%%%%%%%%%%%%%%%%%%%%%%%%%%%%%%%%%%%

\begin{lemma}\label{radiation1}\
Let $c$ be a bounded function on $S^{n-1}$, and $\varphi\in \mathcal{S}_\chi(\R^n)$.
Let $\varepsilon>0$ and $\kappa>0$. 
Suppose that $N_1$ is a positive integer satisfying $N_1\varepsilon\ge \kappa+n-1$.
Then there exists a positive constant $C=C(\varepsilon, \kappa, N_1)$ such that we have
\begin{align*}
& \left|\int_{\Lambda_\varepsilon^1(t, r, \omega)} c\left( %\frac{y-r\theta}{\left|y-r\theta\right|}
t^{-1}(y-r \omega)
\right)\varphi(y) dS^*_y\right|
\\
& \qquad \le C (1+t+r)^{-\kappa}\chi^{-1}(|r-t|)\|c\|_{L^\infty(S^{n-1})}\|\varphi\|_{\chi, 0, N_1}
\end{align*}
for any $(t,r, \omega)\in [0, \infty)\times [0,\infty)\times S^{n-1}$.
\end{lemma}
%%%%%%%%%%%%%%%%%%%%%%%%%%%%%%%%%%%%%%
\begin{proof}
Observing that the total measure of $S(t, r \omega)$ is bounded by $A_n t^{n-1}$,
and that we have
\begin{align*}
|\varphi(y)|\le (1+|t+r|^{2\varepsilon})^{-N_1/2}\chi^{-1}(|r-t|)\|\varphi\|_{\chi, 0,N_1}\\
\le C (1+t+r)^{-\kappa-n+1}\chi^{-1}(|r-t|)\|\varphi\|_{\chi, 0,N_1}
\end{align*}
for any $y\in \Lambda^1_\varepsilon(t, r, \omega)$ with some positive constant $C$,
because $|y|\ge |r-t|$ for any $y\in S(t, r\omega)$.
Thus we obtain the desired result.
\end{proof}
%%%%%%%%%%%%%%%%%%%%%%%%%%%%%%%%%%%%%%%%%%%

%%%%%%%%%%%%%%%%%%%%%%%%%%%%%%%%%%%%%%%%%%%
\begin{lemma}\label{radiation2}\
Let $c\in C^1\left( \overline{B_1} \right)$, 
and $\varphi\in \mathcal{S}_\chi(\R^n)$.
Let $0<\varepsilon\le 1/4$, and $\nu\ge 0$.
Suppose that $N_2$ is a positive integer satisfying $N_2\ge n+2+\nu+(1/\varepsilon)$.
Then there exists a positive constant $C=C(\varepsilon, \nu, N_2)$ such that we have
\begin{align}
\label{radbas02}
& \left|\int_{\Lambda_\varepsilon^2(t, r, \omega)} c\left( %\frac{y-r\theta}{\left|y-r\theta\right|}
t^{-1}(y-r \omega)
\right)\varphi(y) dS^*_y
{}- c(-\omega) {\mathcal R}[\varphi](r-t, \omega)
\right|\\
&\qquad \le C (1+t+r)^{-1}(1+|r-t|)^{-\nu}\chi^{-1}(|r-t|) 
\nonumber\\
& \qquad \qquad \times \|c\|_{C^1\left(\overline{B_1}\right)}\|\varphi\|_{\chi, 1, N_2}
\nonumber
\end{align}
for any $(t, r, \omega)\in[0,\infty)\times [0,\infty)\times S^{n-1}$ 
with $r\ge t/2\ge 1$ and $|r-t|\le (t+r)^\varepsilon$,
where
$$
 \|\psi\|_{C^1\left( \overline{B_1} \right)}=
\sup_{y\in \overline{B_1}} 
\left(|\psi(y)|^2+|\nabla_y \psi(y)|^2
\right)^{1/2}
$$
for $\psi\in C^1\left(\overline{B_1} \right)$.
\end{lemma}
%%%%%%%%%%%%%%%%%%%%%%%%%%%%%%%%%%%%%%%
\begin{proof}
Since the right-hand side of \eqref{radbas02} is invariant under the orthogonal transforms,
we may assume $\omega=e_n$ without loss of generality, where $e_n=(0, \ldots, 0, 1)$.

Suppose $r\ge t/2 \ge 1$, $|r-t|\le (t+r)^\varepsilon$ and
$0<\varepsilon\le 1/4$, in the following.
Then, since $t+r\ge 3$, we get $(t+r)^{\varepsilon-1}\le 3^{-3/4}<1/2$,
which implies $(t+r)^\varepsilon<(t+r)/2$.
If $r> 3t$, then we get $r-t>(t+r)/2>(t+r)^\varepsilon$,
which contradicts the assumption. Hence we obtain
\begin{equation}
 \label{EQU}
 1\le \frac{t}{2}\le r\le 3t,
\end{equation}
and we find that $t$, $r$ and $1+t+r$ are equivalent to each other.

First we prove that
\begin{align}
\label{radbas02+}
& \biggl|\int_{\Lambda_\varepsilon^2(t,r,e_n)} c\left( %\frac{y-r\theta}{\left|y-r\theta\right|}
t^{-1}(y-re_n)
\right)\varphi(y) dS^*_y
\\
& \qquad\qquad\qquad\qquad {}- c(-e_n) 
\int_{\Lambda_\varepsilon^2(t,r,e_n)} \varphi(y) dS^*_y
\biggr| \nonumber\\
& \le 2^{\frac{n-2}2} A_{n-1} 
  \|c\|_{C^1\left(\overline{B_1}\right)}\|\varphi\|_{\chi, 1, N_2}\,t^{-1}
(1+|r-t|^2)^{-\frac{\nu}2}\chi^{-1}(|r-t|).
\nonumber
\end{align}
%for any $(t,r,\theta)\in[0,\infty)\times [0,\infty)\times S^{n-1}$ 
%with $r\ge t/2\ge 1$ and $|r-t|\le (t+r)^\varepsilon$.
We put 
$$
%\lambda_0^\varepsilon=\lambda_0^\varepsilon(t,r)=\frac{(t+r)^{2\varepsilon}-(r-t)^2}{2rt}.
\lambda_0(t,r)=\frac{(t+r)^{2\varepsilon}-(r-t)^2}{2rt}.
$$
Note that $0 \le \lambda_0(t,r) \le 1$.
Writing $t\lambda=y_n-(r-t)$, we find that
$\Lambda_\varepsilon^2(t,r,e_n)$ is equal to
$$
\left\{y=\left(
 t\sqrt{\lambda(2-\lambda)}\zeta, r-t+t\lambda
 \right);
\zeta\in S^{n-2},0\le \lambda \le %\lambda_0^\varepsilon(t,r) 
 \lambda_0(t,r) \right\}.
$$
For the coordinate system $(\zeta, \lambda)$ in the above, we have
\begin{equation}
\label{B1}
dS_y^*=t^{n-1}\lambda^{\frac{n-3}{2}}(2-\lambda)^{\frac{n-3}{2}}d\lambda dS'_\zeta,
\end{equation}
where $dS'_\zeta$ denotes the area element on $S^{n-2}$.
We also note that
\begin{equation}
\label{B2}
|y|^2=(r-t)^2+2rt\lambda.
\end{equation}
We put $\eta(y, r)=|y-re_n|^{-1}(y-r e_n)\in S^{n-1}$.
Then we get
\begin{equation*}
|\eta(y,r)-(-e_n)|=\sqrt{2+2\eta_n(y,r)}=\sqrt{2\lambda}
\end{equation*}
for any $y\in \Lambda_\varepsilon^2(t,r,e_n)$.
Hence, by the mean value theorem, we get
\begin{align}
%\label{C1}
%& 
\left| \int_{\Lambda_\varepsilon^2(t,r, e_n)} \left\{c\left(\eta(y,r)\right)-c(-e_n)\right\} \varphi(y)dS_y^*
\right|
%\\ & 
\quad \le \sqrt{2} \|c\|_{C^1\left( \overline{B_1}\right)} J(t,r), \nonumber
\end{align}
where we put
$$
J(t,r)=\int_{\Lambda_\varepsilon^2(t,r, e_n)}\lambda^{\frac{1}{2}} |\varphi(y)| dS_y^*
$$
with $t\lambda=y_n-(r-t)$.
Recalling \eqref{B1} and \eqref{B2}, we obtain
\begin{align}
\label{J1}
& \chi(|r-t|)J(t,r)\\
& \quad \le A_{n-1} t^{n-1} \|\varphi\|_{\chi, 0, N_2} 
\int_0^{\lambda_0(t,r)}
\frac{\lambda^{\frac{n-2}{2}}(2-\lambda)^{\frac{n-3}{2}}}{(1+(r-t)^2+2rt\lambda)^{\frac{N_2}{2}}}
d\lambda
\nonumber\\
& \quad \le 2^{-\frac{1}{2}}A_{n-1} r\left(\frac{t}{r}\right)^{\frac{n}{2}} \|\varphi\|_{\chi, 0,N_2}
 \nonumber\\
& \qquad\quad \times \int_0^{\lambda_0(t,r)}
\frac{(2rt\lambda)^{\frac{n-2}{2}}}{(1+|r-t|^2)^{\frac{\nu}{2}}(1+2rt\lambda)^{\frac{N_2-\nu}{2}}}d\lambda
\nonumber\\
& \quad \le 2^{\frac{n-1}{2}}A_{n-1}\,r(1+|r-t|^2)^{-\frac{\nu}{2}} \|\varphi\|_{\chi, 0,N_2}
\int_0^{\infty}\frac{1}{(1+2rt\lambda)^{2}}d\lambda
\nonumber\\
& \quad \le 2^{\frac{n-3}{2}}A_{n-1}\,t^{-1}(1+|r-t|^2)^{-\frac{\nu}{2}} \|\varphi\|_{\chi, 0,N_2}.
\nonumber
\end{align}
This estimate yields \eqref{radbas02+} immediately.

By \eqref{radbas02+}, we find that, 
in order to show \eqref{radbas02}, it suffices to prove
\begin{align}
\label{C3}
& \left |\int_{\Lambda_\varepsilon^2(t,r, e_n)} \varphi(y) dS_y^*
    -{\mathcal R}[\varphi](r-t, e_n) \right|
\\
&\quad \le C\,t^{-1}(1+|r-t|)^{-\nu}\chi^{-1}(|r-t|)\|\varphi\|_{\chi, 1, N_2}
\nonumber
\end{align}
with some positive constant $C$.
We observe that
\begin{align}
\left| \int_{\Lambda_\varepsilon^2(t,r, e_n)} \lambda \varphi(y) dS_y^*\right|
\le C\,(rt)^{-1}(1+|r-t|)^{-\nu}\chi^{-1}(|r-t|)\|\varphi\|_{\chi, 0, N_2},
\nonumber
\end{align}
which can be shown similarly to \eqref{J1}.
%Similarly to  we see that  $|I_0(t,r)-I_1(t,r)|$ is bounded from above by
%$$
%2^{\frac{n-3}{2}}A_{n-1} \|\varphi\|_{\chi, 0, N_2}(1+|r-t|^2)^{-\frac{\nu}{2}} 
%\chi^{-1}(|r-t|)\int_0^\infty 
%\frac{1}{(1+2rt\lambda)^{\frac{3}{2}}} d\lambda,
%$$
%whose last integral is bounded by $(rt)^{-1}$.
Therefore, \eqref{C3} follows from
\begin{align}
\label{C3+}
& |I_1(r,t)-{\mathcal R}[\varphi](r-t, e_n)|
\\
&\quad \le C\,t^{-1}(1+|r-t|)^{-\nu}\chi^{-1}(|r-t|)\|\varphi\|_{\chi, 1, N_2},
\nonumber
\end{align}
where we put
$$
I_1(r,t)=\int_{\Lambda_\varepsilon^2(t,r, e_n)} (1-\lambda) \varphi(y) dS_y^*.
$$

Introducing a new coordinate $\rho=t\sqrt{\lambda(2-\lambda)}$,
we get
\begin{align*}
I_1(t,r)
=& \int_{S^{n-2}} \left(
\int_0^{\rho_0(t,r)}
 \varphi \left(\rho \zeta, r-\sqrt{t^2-\rho^2} \right) 
\rho^{n-2}d\rho \right) dS'_\zeta,
\end{align*}
where 
$ %\rho_0^\varepsilon(t,r)=
   %t\sqrt{\lambda_0^\varepsilon(t,r)(2-\lambda_0^\varepsilon(t,r))}
   \rho_0(t,r)=
        t\sqrt{\lambda_0(t,r)(2-\lambda_0(t,r))}$.
While, we have
\begin{align}
\label{C3++}
{\mathcal R}[\varphi](r-t, e_n)= %\int_{y\cdot e_n=r-t} \varphi(y)dS_y=
I_2(t,r)+I_3(t,r),
\end{align}
where we put
\begin{align*}
& I_2(t,r)=\int_{S^{n-2}}\left( %\int_0^{\rho_0^\varepsilon(t,r)} 
\int_0^{\rho_0(t,r)}
\varphi\left(\rho\zeta, r-t\right)\rho^{n-2} d\rho \right) dS_\zeta',
\\
& I_3(t,r)=\int_{S^{n-2}}\left( %\int_{\rho_0^\varepsilon(t,r)}^\infty 
\int_{\rho_0(t,r)}^\infty
\varphi(\rho\zeta, r-t)\rho^{n-2} d\rho\right)dS_\zeta'.
\end{align*}
Since $t-\sqrt{t^2-\rho^2}=
\rho^2\left(t+\sqrt{t^2-\rho^2}\right)^{-1}$, we get
\begin{align*}
& \left|
\varphi\left(\rho\zeta, r-\sqrt{t^2-\rho^2}\right)
{}-\varphi(\rho\zeta, r-t)
\right|\\
& \quad 
\le \left(t-\sqrt{t^2-\rho^2}\right) \int_0^1 \left|(\pa_n \varphi)
\left(\rho\zeta, r-t+\tau\left(t-\sqrt{t^2-\rho^2}\right)\right)\right|d\tau 
%\\
%& \quad 
%\le \frac{\|\varphi\|_{\chi, 1, N_2} \rho^2}{t(1+|r-t|^2)^{\frac{\nu}{2}}(1+\rho^2)^{\frac{N_2-\nu}{2}} \chi(|r-t|)}
\\
& \quad 
\le \frac{\rho^2 \|\varphi\|_{\chi, 1,N_2}}{t(1+|r-t|^2)^{\frac{\nu}{2}}(1+\rho^2)^{\frac{n+2}{2}}\chi(|r-t|)}
\end{align*}
for $0\le \rho\le \rho_0(t,r)$,
which yields
\begin{align}
\label{d1}
& |I_1(t,r)-I_2(t,r)|\\
& \qquad 
\le 
\frac{A_{n-1} \|\varphi\|_{\chi, 1,N_2}}{t(1+|r-t|^2)^{\frac{\nu}2}\chi(|r-t|)}
\int_0^\infty \frac{\rho}{(1+\rho^2)^{\frac{3}{2}}}
d\rho \nonumber\\
% CA_{n-1}\|\varphi\|_{\chi, 1,N_2}t^{-1}(1+|r-t|)^{-\nu}\chi^{-1}(|r-t|)
% \nonumber\\
% & \qquad\qquad \times\int_0^\infty \frac{\rho}{(1+\rho^2)^{\frac{3}{2}}}
% d\rho \nonumber\\
& \qquad 
\le CA_{n-1} \|\varphi\|_{\chi, 1,N_2}\,t^{-1}(1+|r-t|)^{-\nu}\chi^{-1}(|r-t|),
\nonumber
\end{align}
where $C$ is a positive constant depending only on $\nu$.

Finally, we evaluate $I_3(t, r)$. Notice that
\begin{equation}
\label{D2}
\left(\rho_0(t,r)\right)^2+(r-t)^2\ge \frac{3}{32}(t+r)^{2\varepsilon}.
\end{equation}
In fact, it is trivial when $(t+r)^{2\varepsilon}/2\le (r-t)^2
\left(\le (t+r)^{2\varepsilon}\right)$. On the other hand,
when $(r-t)^2\le (t+r)^{2\varepsilon}/2$, recalling that we have $(t+r)^\varepsilon<(t+r)/2$,
we get
\begin{align*}
\left(\rho_0(t,r)\right)^2=& \frac{(t+r)^2-(t+r)^{2\varepsilon}}{4r^2}
\left\{(t+r)^{2\varepsilon}-(r-t)^{2}\right\}\ge\frac{3}{32}(t+r)^{2\varepsilon},
\end{align*}
which shows \eqref{D2}. We thus find %By \eqref{D2}, we conclude
\begin{align}
\label{I3}
|I_3(t,r)|\le & \frac{A_{n-1} \|\varphi\|_{\chi, 0, N_2}}{\chi(|r-t|)} 
\int_{\rho_0(t,r)}^\infty
\frac{\rho}{(1+\rho^2+(r-t)^2)^{\frac{\nu+5+(1/\varepsilon)}{2}}} d\rho \\
\le & C\|\varphi\|_{\chi, 0, N_2} (t+r)^{-1}(1+|r-t|)^{-\nu}\chi^{-1}(|r-t|), 
\nonumber
\end{align}
where $C$ is a constant depending only on $\nu$, $n$ and $\varepsilon$.
Now, \eqref{C3+} follows from \eqref{C3++}, \eqref{d1} and \eqref{I3}.
This completes the proof.
\end{proof}

\section{Proof of Theorem \ref{radiation}}\label{ProofRadiation}
%Finally we are going to give the proof of Theorem~\ref{radiation}.\\
%{\it Proof of Theorem~\ref{radiation}.} 
To begin with, we note that, for any $A\in \R$ and any nonnegative integer $k$,
there exists a positive constant $C$ such that we have
$$
C^{-1} e^{A\jb{x}} \sum_{|\alpha|\le k} |\pa_x^\alpha \psi(x)|
\le \sum_{|\alpha|\le k} \left|\pa_x^\alpha \bigl(e^{A\jb{x}} \psi(x)\bigr)\right|
\le C e^{A\jb{x}} \sum_{|\alpha|\le k} |\pa_x^\alpha \psi(x)|
$$
for any $x\in \R^n$ and any $\psi\in C^\infty(\R^n)$.
In fact, the latter half is almost apparent, 
and the first half is nothing but the latter half with $A$ and $\psi(x)$
being replaced by $-A$ and $e^{A\jb{x}}\psi(x)$, respectively.
% and the latter half implies
% $$
% \sum_{|\alpha|\le k} \left|\pa_x^\alpha \bigl(e^{-A\jb{x}}e^{A\jb{x}}\psi(x)\bigr)\right|\le 
% C e^{-A\jb{x}}\sum_{|\alpha|\le k} \left|\pa_x^\alpha \bigl(e^{A\jb{x}} \psi(x)\bigr)\right|,
% $$
% which immediately leads to the first half.

Let the assumptions in Theorem~\ref{radiation}
be fulfilled.
Then, by Theorem~\ref{inverse}, there exists $f_+\in \Hi^\infty(\R^n)$ satisfying \eqref{iw0} and \eqref{iw01}.
We write $u(t,\cdot)$ and $u_+(t, \cdot)$ for the first components
of $U(t)\vec{f}$ and $U_0(t)\vec{f}_+$, respectively.

First we claim that we have
\begin{align}
\label{rp01}
\left\|\exp(3\mu\jb{\cdot}/4) \left(\Gamma^\alpha u(t,\cdot)-\Gamma^\alpha u_+(t,\cdot)\right) \right\|_{H^{[n/2]+1}(\Omega)} \qquad\qquad
\\
\le C\exp(-\mu t)\|\vec{f}\|_{{\mathcal
H}^{[n/2]+|\alpha|}(\Omega)},\quad t\ge 0
\nonumber
\end{align}
for any multi-index $\alpha$.
Let $\Gamma^\alpha=\pa_t^j\pa_x^\beta O^\gamma$ with a
nonnegative integer $j$, and multi-indices $\beta$, $\gamma$,
and let $\alpha^\prime$ be a multi-index satisfying $|\alpha'| \le [n/2]+1$.
If $j$ is even, then we get 
\begin{align*}
 |\pa_x^{\alpha'}\Gamma^\alpha (u-u_+)(t,x)|=& |\Delta^{j/2}\pa_x^{\alpha'+\beta} O^\gamma
  (u-u_+)(t,x)|\\
 \le & C(1+|x|)^{|\gamma|} \sum_{|\beta'|\le |\alpha|+|\alpha'|}
 |\pa_x^{\beta'}(u-u_+)(t,x) |\\
 \le & Ce^{\mu \jb{x}/4}\sum_{|\beta'|\le |\alpha|+|\alpha'|}
 |\pa_x^{\beta'}  (u-u_+)(t,x)|
\end{align*}
with some positive constant $C=C(\alpha, \alpha', \mu)$.
Similarly, if $j$ is odd, we get
$$
 |\pa_x^{\alpha'} \Gamma^\alpha (u-u_+)(t,x)|
\le C e^{\mu \jb{x}/4} 
\sum_{|\beta'|\le |\alpha|+|\alpha'|-1}
 |\pa_x^{\beta'}\pa_t(u-u_+)(t,x)|.
$$
Hence the left-hand side of \eqref{rp01} is bounded by
$$
C\|\exp(\mu\jb{\cdot})(U(t)\vec{f}-U_0(t)\vec{f}_+)\|_{{\mathcal H}^{[n/2]+|\alpha|}(\Omega)},
$$
and \eqref{iw0} implies \eqref{rp01}.

By \eqref{rp01} and the Sobolev imbedding theorem, we get
% \begin{align}
% & |\pa_x^\alpha (U(t)\vec{f}-U_0(t)\vec{f}_+)(x)|\\
% & \quad \le C_\alpha e^{-\mu \jb{x}}\sum_{|\beta|\le |\alpha|} |\pa_x^\beta(e^{\mu\jb{x}}
% (U(t)\vec{f}-U_0(t)\vec{f}_+)(x)| \nonumber\\
% & \quad \le CC_\alpha e^{-\mu (t+\jb{x})}\|\vec{f}\|_{\Hi^{k+[n/2]+1}(\Omega)},
% \quad  |x|\ge t/2 \ge 1\nonumber
% \end{align}
% for $|\alpha|\le k$, where $C_\alpha$ is a positive constant depending on $\alpha$.
% Hence we obtain
\begin{align}
\label{Finale01}
& 
\left|\Gamma^\alpha (u-u_+)(t,x) \right|
\\
  & \quad \le C_k e^{-3\mu \jb{x}/4}
%\sum_{|\beta|\le |\alpha|} 
 \left\|e^{3\mu\jb{\cdot}/4}\Gamma^\alpha
 (u-u_+)(t, \cdot)\right\|_{H^{[n/2]+1}(\Omega)} \nonumber\\
& \quad 
\le C_ke^{-3\mu(t+\jb{x})/4}\|\vec{f}\|_{\Hi^{k+[n/2]}(\Omega)}
\nonumber\\
& \quad 
\le C_k (1+t+|x|)^{-(n+1)/2}e^{-\mu(t+|x|)/2} \|\vec{f}\|_{\Hi^{k+[n/2]}(\Omega)}
\nonumber
\end{align}
for $|x|\ge t/2\ge 1$ and $|\alpha|\le k$, 
where $C_k$ is a positive constant.
Hence we find that our task is to show \eqref{ki6} and \eqref{ki7} with
$u$ being replaced by $u_+$.

From \eqref{iw01} we see that 
$\vec{f}_+\in \left({\mathcal S}_\chi(\R^n)\right)^2$ with
$\chi(s)=\exp(\mu s/2)$.
Therefore Proposition~\ref{radiation0} with $\nu=0$ immediately implies
\eqref{ki6} and \eqref{ki7} with
$u$ being replaced by $u_+$.
This completes the proof. \qed

% \noindent
% {\bf Concluding Remark}.
% Proposition~\ref{radiation0} itself, which is the general result for rapidly decreasing data,
% seems  new and worth stating explicitly, and thus
% we give Corollary~\ref{radiation} in the present form. However, since \eqref{iw01}
% shows that $\vec{f}$
% has a far better decay property than we can expect in general for rapidly decreasing functions,
% we can improve the result of Corollary~\ref{radiation} in the following way.
% By apparent modification of the proof, we can show that
% if $\varphi \in C^\infty(\R^n)$ satisfies
% $$
% \sup_{x\in \R^n} e^{\mu |x|} \sum_{|\alpha|\le k} \left|\pa_x^\alpha \varphi (x)\right|<\infty
% $$
% for any nonnegative integer $k$ with some $\mu>0$ being independent of $k$,
% then, instead of \eqref{DecayRadon},  we have
% $$
% \sup_{(s, \omega)\in \R\times S^{n-1}} e^{\mu |s|/2} |\pa_s^jo^\alpha {\mathcal R}[\varphi](s,\omega)|<\infty 
% $$
% for any nonnegative integer $j$ and any multi-index $\alpha$.
% We can also show that if $\vec{f}\in \left(C^\infty(\R^n)\right)^2$ satisfies
% $\sup_{x\in \R^n} e^{\mu|x|} |\vec{f}(x)|_k<\infty$
% for any nonnegative integer $k$ with some $\mu>0$,
% then \eqref{ki60} and \eqref{ki70} with $(1+|r-t|)^{-\nu}$ 
% being replaced by $e^{-\mu |r-t|/2}$ are valid. 
% Here $|\vec{f}|_k$ is given by \eqref{Not01}.
% Accordingly, in view of  \eqref{iw0} and \eqref{iw01}, we find that Corollary~\ref{radiation} with
% $(1+|r-t|)^{-\nu}$ in \eqref{ki6} and \eqref{ki7} being replaced by
% $e^{-\mu |r-t|/2}$ is true.

\end{document}